\documentclass{article}
\usepackage{arxiv}
\usepackage{amsfonts,amssymb,amsmath,amsthm}
\usepackage{graphicx}
\usepackage{algorithm,algorithmic,multirow}
\usepackage{epstopdf}
\usepackage[labelformat=simple]{subfig}

\usepackage{xcolor}
\usepackage{latexsym,bm}        % for things like \Box and \Diamond
\graphicspath{ {./figs/} }

\newtheorem{theorem}{Theorem}[section]
\newtheorem{remark}{Remark}[section]
\DeclareMathOperator{\sech}{sech}

\usepackage{hyperref}
\makeindex

\title{Global energy preserving model reduction for multi-symplectic PDEs}

\author{
  \href{https://orcid.org/0000-0001-5262-063X}{\includegraphics[scale=0.06]{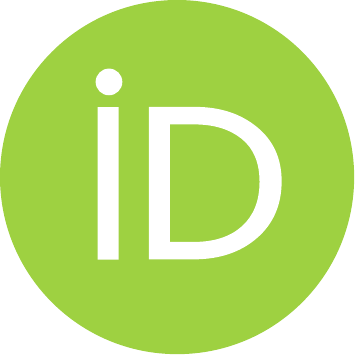}\hspace{1mm}Murat Uzunca} \\
   Department of Mathematics\\
	Sinop University, Sinop-Turkey\\
     \texttt{muzunca@sinop.edu.tr}\\
		    \And
     \href{https://orcid.org/0000-0003-1037-5431}{\includegraphics[scale=0.06]{orcid.pdf}\hspace{1mm}B\"ulent Karas\"ozen} \\
     Institute  of Applied Mathematics \& Department of Mathematics\\
     Middle East Technical University, Ankara-Turkey\\
     \texttt{bulent@metu.edu.tr}
          \And
\href{https://orcid.org/0000-0002-0837-9364}{\includegraphics[scale=0.06]{orcid.pdf}\hspace{1mm}Ayhan Ayd{\i}n} \\
Department of Mathematics\\
At{\i}l{\i}m  University, Ankara-Turkey\\
	\texttt{ayhan.aydin@atilim.edu.tr}
}

\begin{document}

\maketitle

\begin{abstract}
Many Hamiltonian systems can be recast in multi-symplectic form. We develop a reduced-order model (ROM) for multi-symplectic Hamiltonian partial differential equations (PDEs) that preserves the global energy. The full-order solutions are obtained by finite difference discretization in space and the global energy preserving average vector field (AVF) method. The ROM is constructed in the same way as the full-order model (FOM) applying proper orthogonal decomposition (POD) with the Galerkin projection. The reduced-order system has the same structure as the FOM, and preserves the discrete reduced global energy. Applying the discrete empirical interpolation method (DEIM), the reduced-order solutions are computed efficiently in the online stage. A priori error bound is derived for the DEIM approximation to the nonlinear Hamiltonian. The accuracy and computational efficiency of the ROMs are demonstrated for the Korteweg de Vries (KdV) equation, Zakharov-Kuznetzov (ZK) equation, and nonlinear Schr{\"o}dinger (NLS) equation in multi-symplectic form. Preservation of the reduced energies shows that the reduced-order solutions ensure the long-term stability of the solutions.
\end{abstract}


%%%%%%%%%%%%%%%%%%%%%%%%%%%%%%%%%%%%%%%%%%%%%%%%%%%%%%%%%%%%%%%%%%%%%%%%%%%%%%%%%%%%%%%%%%%%
%%%%%%%%%%%%%%%%%%%%%%%%%%%%%%%%%%%%%%%%%%%%%%%%%%%%%%%%%%%%%%%%%%%%%%%%%%%%%%%%%%%%%%%%%%%%
\section{Introduction}

Many partial differential equations (PDEs) can be recast in multi-symplectic Hamiltonian formulation such as the Sine-Gordon equation \cite{Reich04}, the Korteweg-de Vries (KdV) equation \cite{Cai19,Eidnes20,Gong14}, the Camassa-Holm equation \cite{Cohen08ch}, the nonlinear Schr{\"o}dinger (NLS) equation \cite{Cai19,Gong14}, the Maxwell's equations \cite{Cai16}, the Boussinesq equation \cite{Aydin08} and so on.
The solution of multi-symplectic PDEs admits a multi-symplectic conservation
law (MSCL), local/global energy and momentum conservation laws.
Conservation laws play an important role in conservative PDEs, especially in the theory of solitons.
Numerical integrators that preserve the geometric structure of dynamical systems such as symplecticity and multi-symplecticity, invariants such as the energy and momentum lead to stable solutions in long-term \cite{Hairer10gni, Reich05}.
A numerical method that satisfies a discrete version of MSCL law is called multi-symplectic integrator \cite{Bridges01a}.
In the last two decades, many multi-symplectic integrators have been developed; the box/Preissmann scheme \cite{Reich04}, the Euler-box scheme \cite{Moore03}, the Fourier pseudo-spectral collocation scheme \cite{Bridges01a}, wavelet collocation method \cite{Zhu10}.
Local and global energy/momentum preserving methods have been developed in \cite{Cai19,Gong14} with the average vector field (AVF) method and in \cite{Eidnes20,Cai20} with the linearly implicit Kahan's method.

Real-time simulations of PDEs require a large amount of computer memory and computing time
with fully resolved simulations multiple times for different parameter settings or in multi-query scenarios such as in optimization.
In the last two decades, many model order reduction (MOR) methods have been developed that allow the construction of low-dimensional reduced-order models (ROMs) for the high fidelity full-order models (FOMs).
The solutions of the high fidelity FOMs are projected on low dimensional reduced spaces usually using the proper orthogonal decomposition (POD) method \cite{Berkooz93,Sirovich87}, which is a widely used ROM technique.
Conservation of physical quantities like energy is not, in general, guaranteed by the ROMs.
When the physical quantities are violated, ROMs frequently result in an unstable or qualitatively wrong reduced system, even when the high-fidelity system is stable.
The stability of reduced models over long-time integration has been investigated in the context of Lagrangian systems \cite{Carlberg15}, and for port-Hamiltonian systems, \cite{Chaturantabu16}.
For canonical Hamiltonian systems, like the linear wave equation, Sine-Gordon equation, NLS equation, singular value decomposition (SVD) based symplectic model reduction techniques with POD Galerkin projection are constructed by orthogonal bases \cite{Peng16} and non-orthogonal bases \cite{Buchfink18} that capture the symplectic structure of Hamiltonian systems to ensure long term stability of the reduced model. For parametric Hamiltonian systems, symplectic bases are generated using greedy approaches in \cite{Hesthaven16,Buchfink20greedy}.
Parallel to these, energy (Hamiltonian) preserving ROMs have been developed for non-canonical Hamiltonian PDEs like the KdV equation \cite{Gong17,Miyatake19,Karasozen21}, the NLS equation \cite{Karasozen18nls}, and with state-dependent Poisson structure such as  rotating shallow water equations \cite{Karasozen21sw,Karasozen22rtswe}. These are all global ROM techniques that  maintain the globalized properties such as the symplectic structure or the Hamiltonian of the full order data in the reduced order representation. On the other hand, because of the nondissipative phenomena, the Hamiltonian PDEs do not possess a global low-rank structure, i.e., they require a sufficiently large approximation space to achieve
accurate solutions when a global ROM is used.
This can be ascribed to the local low-rank nature of nondissipative phenomena, like advection
and wave-type problems, characterized by slowly decaying Kolmogorov $n$-widths for transport equations.
Recently, localized ROM have been developed for Hamiltonian PDEs.
A reduced basis method is developed in \cite{Hesthaven21} for non-canonical Hamiltonian systems that preserves general Poisson structure by "freezing"  the phase space manifold structure in each discrete temporal interval,  then recasts the local problem in canonical
form. In this way, a local reduced model is constructed in canonical Hamiltonian
form. In \cite{Hesthaven22}, a rank-adaptive structure-preserving dynamical reduced-based method is constructed based on a residual error estimator. The FOMs are approximated
on local reduced spaces that are adapted in time using dynamical low-rank approximation
techniques. In \cite{Pagliantini21}, nonlinear structure-preserving model reduction is proposed where the reduced phase space evolves in time.
The reduced system is obtained by a symplectic projection of the Hamiltonian vector field onto the tangent space of the approximation manifold at each reduced state as in the dynamical low-rank approximation.
We refer to \cite{Hesthaven22an} for an overview about the structure-preserving ROMs for Hamiltonian systems.

In this paper, a global ROM is investigated for multi-symplectic PDEs that preserve the energy/Hamiltonian globally. To the best of our knowledge, ROMs for the multi-symplectic PDEs are not considered before.
The FOM is constructed by discretizing the PDE in multi-symplectic form by finite differences in space and integrating in time with global energy preserving AVF method \cite{Gong14,Cai20}. The reduced system obtained with the POD has the same multi-symplectic structure as the FOM with a reduced energy/Hamiltonian. Following the ROM technique
in \cite{Gong14, Miyatake19} and applying Galerkin projection, the reduced ordinary differential equations (ODEs) are integrated in time again with the AVF method.
In this way, the discrete reduced global energy is preserved with the POD-ROM (P-ROM) exactly. An important feature of the ROMs is the offline-online decomposition of computations.
The computation of the FOM and the construction of the reduced basis are performed in the offline stage, whereas the reduced system is solved by projecting the problem onto the low-dimensional reduced space in the online stage.
To accelerate the computation of the nonlinear terms in the reduced form, hyper-reduction techniques such as the discrete empirical interpolation method (DEIM) \cite{chaturantabut10nmr,Gugercin16a}, are used, i.e., the reduced-order solutions do not depend on the FOM, and an efficient offline-online decomposition is achieved.
An a priori error bound is derived for the POD-DEIM approximated ROM, which depends on the number of DEIM modes. Numerical tests on one- and two-dimensional multi-symplectic Hamiltonian PDEs such as KdV equation, NLS equation, Zakharov-Kuznetzov (ZK) equation, and two-dimensional NLS equation with soliton solutions, show that the ROMs have the same level of accuracy. The discrete reduced global energies are accurately preserved with an increasing number of POD and DEIM modes with small oscillation over time, which ensures the long-term stability of the solutions. Speedups of the computing times of the ROMs over the FOMs of moderate order are achieved depending on the problem type, which is typical for problems with wave and soliton solutions.

The global energy preserving ROM investigated in this paper does not possess a global low-rank structure as the existing structure preserving ROM techniques for Hamiltonian systems. On the other side, the numerical tests on one- and two-dimensional multi-symplectic Hamiltonian PDEs show that the structure-preserving approach proposed in
this work requires a moderately large number of POD/DEIM modes to achieve robust and efficient reduced solutions. Localized ROM approaches can be extended for the multi-symplectic Hamiltonian PDEs, by which the local low-rank structure of the wave and soliton type solution can be preserved using a small number of POD/DEIM modes, but this might require more computing time.
Similar to the global energy preserving ROMs, global momentum preserving ROMs can also be developed from the momentum preserving FOMs as in \cite{Gong14}.
The multi-symplectic PDEs also possess local conservation laws, such as local multi-symplectic structure, energy, and momentum conservation. Therefore, it would be better suited to develop ROM techniques to preserve the local properties of the multi-symplectic PDEs in the reduced-order form, which will be the subject of our future research.

The paper is organized as follows. In Section~\ref{sec:mcpde}, we introduce general multi-symplectic PDEs with their underlying properties.  The global energy preserving FOM in space and time is introduced in Section~\ref{sec:fom}. The global energy preserving ROMs with the POD and the DEIM are constructed in Section~\ref{sec:rom}. In Section~\ref{sec:num}, results of numerical tests are presented and compared with those in the literature.  Some conclusions are drawn in Section~\ref{sec:conc}.

%%%%%%%%%%%%%%%%%%%%%%%%%%%%%%%%%%%%%%%%%%%%%%%%%%%%%%%%%%%%%%%%%%%%%%%%%%%%%%%%%%%%%%%%%%%%
%%%%%%%%%%%%%%%%%%%%%%%%%%%%%%%%%%%%%%%%%%%%%%%%%%%%%%%%%%%%%%%%%%%%%%%%%%%%%%%%%%%%%%%%%%%%
\section{Multi-symplectic PDEs}
\label{sec:mcpde}

The multi-symplectic formulation of a Hamiltonian PDE is given as \cite{Cai19,Bridges01a}
\begin{align}\label{ms1}
K\partial_tz + \sum_{s=1}^d L^s\partial_{x_s}z = \nabla _z S(z), \quad (\bm{x},t)\in \Omega\times (0,T],
\end{align}
where $z(\bm{x},t)=(z_{(1)}(\bm{x},t),\ldots ,z_{(m)}(\bm{x},t))^T$ is the $m$-component  state vector, $d$ denotes the space dimension, $K$ and $L^s$ are skew-symmetric matrices of size $m\times m$, $s=1,\ldots,d$, $\bm{x}=(x_1,\ldots ,x_d)^T\in\Omega\subset\mathbb{R}^{d}$, and $S(z):\mathbb{R}^{m}\mapsto \mathbb{R}$ is a smooth functional.
In \eqref{ms1}, $\partial_t$ and $\partial_{x_s}$ stand for the partial differentiation with respect to the time $t$ and the $s$th space coordinate $x_s$, respectively.

The Hamiltonian system \eqref{ms1} admits MSCL \cite{Cai19,Gong14,Bridges01a}
\begin{equation}\label{mscl}
\partial_t \omega + \sum_{s=1}^d\partial_{x_s}\kappa_s =0,
\end{equation}
with the differential $2$-forms given by
$$
\omega = dz\wedge K_+dz \; , \qquad \kappa_s = dz\wedge L^s_+dz, \quad s=1,\ldots ,d,
$$
and local energy conservation law (LCL)
\begin{equation}\label{lcl}
\partial_tE + \sum_{s=1}^d\partial_{x_s}F^s = 0, \qquad E(z) = S(z) + \sum_{s=1}^d (\partial_{x_s}z)^TL^s_+z, \quad F^s = -(\partial_tz)^TL^s_+z,
\end{equation}
where the matrices $K_+$ and $L^s_+$ are
$$
K = K_+ - K_+^T, \quad L^s = L^s_+-(L^s_+)^T.
$$

Under  periodic and homogeneous Dirichlet boundary conditions, the proposed LCL \eqref{lcl} leads to the global energy conservation law (GCL)
\begin{equation}\label{energy}
\partial_t\varepsilon =0, \qquad \varepsilon (t)=\int_{\Omega} E(z(t))d\bm{x}.
\end{equation}

In the last two decades, many geometric integrators have been proposed  for the multi-symplectic PDEs  that preserve the MSCL \eqref{mscl}, LCL \eqref{lcl} and/or GCL \eqref{energy} \cite{Reich04,Cai19,Eidnes20,Gong14,Bridges01a,Moore03,Zhu10,Cai20}. In Section~\ref{sec:fom},  we consider a FOM that preserves the discrete form of the global energy \eqref{energy} under periodic boundary conditions.

Throughout the paper, we consider two well-known multi-symplectic PDEs; the KdV equation and the NLS equation.
The KdV equation is given by
\begin{equation}\label{kdv1}
u_t + \eta u u_x+\gamma^2 u_{xxx}  = 0,
\end{equation}
where $\eta, \gamma \in \mathbb{R}$ are real parameters. Introducing the potential $\phi_x=u$, momenta $v=\gamma u_x$, and $w=\phi_t/2 + \gamma v_x + \eta u^2/2$, yields
\begin{equation}\label{kdvms1}
\begin{aligned}
\frac{1}{2}u_t+w_x&=0,\\
-\frac{1}{2}\phi_t-\gamma v_x&=-w+\frac{\eta}{2}u^2,\\
\gamma u_x&=v,\\
-\phi_x&=-u.
\end{aligned}
\end{equation}
Then, the multi-symplectic formulation \eqref{ms1} for the KdV equation \eqref{kdv1} is given  with $z = (\phi,u,v,w)^T$, the Hamiltonian
$S(z)=v^2/2-uw+\eta u^3/6$, and the skew-symmetric matrices
\begin{equation*}
K=\begin{bmatrix}
    0       &\frac{1}{2} & 0 & 0 \\
    -\frac{1}{2}       & 0 & 0 & 0 \\
      0      & 0 & 0 & 0 \\
      0      & 0 & 0 & 0  \\
\end{bmatrix},\qquad
L=
\begin{bmatrix}
    0       &0 & 0 &1  \\
   0       & 0 & -\gamma & 0  \\
      0      & \gamma & 0 & 0  \\
      -1      & 0 & 0 & 0 \\
\end{bmatrix}.
\end{equation*}
Under periodic boundary conditions, the system \eqref{kdvms1} possesses the global energy conservation law
\begin{equation}\label{kdvenergy1}
\varepsilon (t) = \int_{\Omega} \left(  \frac{\eta}{6}u^3 - \frac{\gamma^2}{2}(u_x)^2 \right)dx.
\end{equation}

The one-dimensional NLS equation is given by
\begin{equation}\label{nls1}
iu_t + u_{xx} + \beta |u|^2u  = 0,
\end{equation}
with the complex-valued  solution $u(x,t)$, and a positive parameter $\beta>0$. For real-valued functions $p(x,t)$ and $q(x,t)$, decomposing the solution as $u=p+iq$, and introducing the momenta $v=p_x$ and $w=q_x$, we get
\begin{equation}\label{nlsms1}
\begin{aligned}
q_t-v_x &= \beta (p^2+q^2)p,\\
-p_t-w_x &= \beta (p^2+q^2)q,\\
p_x&=v,\\
q_x&=w.
\end{aligned}
\end{equation}
Then, the multi-symplectic formulation \eqref{ms1} for the NLS equation \eqref{nls1} is given with $z = (p,q,v,w)^T$, the Hamiltonian
$S(z)=(v^2+w^2)/2 + \beta (p^2+q^2)^2/4$, and the skew-symmetric matrices
\begin{equation*}
K=\begin{bmatrix}
    0       & 1 & 0 & 0 \\
    -1  & 0 & 0 & 0 \\
      0      & 0 & 0 & 0 \\
      0      & 0 & 0 & 0  \\
\end{bmatrix},\qquad
L=
\begin{bmatrix}
   0   & 0 & -1 & 0  \\
   0   & 0 & 0  & -1  \\
   1   & 0 & 0  & 0  \\
   0   & 1 & 0  & 0 \\
\end{bmatrix}.
\end{equation*}
Under periodic boundary conditions, the global energy
\begin{equation}\label{nlsenergy1}
\varepsilon (t) = \frac{1}{2}\int_{\Omega} \left(  \frac{\beta}{2}(p^2+q^2)^2 - (p_x)^2  - (q_x)^2 \right)dx,
\end{equation}
is preserved.

%%%%%%%%%%%%%%%%%%%%%%%%%%%%%%%%%%%%%%%%%%%%%%%%%%%%%%%%%%%%%%%%%%%%%%%%%%%%%%%%%%%%%%%%%%%%
%%%%%%%%%%%%%%%%%%%%%%%%%%%%%%%%%%%%%%%%%%%%%%%%%%%%%%%%%%%%%%%%%%%%%%%%%%%%%%%%%%%%%%%%%%%%
\section{Full-order model}
\label{sec:fom}

The FOM is constructed by discretizing the multi-symplectic PDE \eqref{ms1} by finite differences in space and integrating in time with the global energy preserving AVF method \cite{Gong14,Cai20}.
In this section, we give the discrete formulations for the following one-dimensional ($d=1$) multi-symplectic PDE
\begin{align}\label{ms}
K\partial_tz + L\partial_{x}z &= \nabla _z S(z), & (x,t)\in [x_L,x_R]\times (0,T],
\end{align}
where $z(x,t)=(z_{(1)}(x,t),\ldots ,z_{(m)}(x,t))^T$ is the vector of state variables.

We consider uniform partition of the spatial domain $[x_L,x_R]$ as $x_j=x_L+(j-1)\Delta x$, $j=1,\ldots ,N+1$, with the mesh size $\Delta x=(x_R-x_L)/N$ and of the time domain $[0,T]$ as  $t_k=k\Delta t$, $k=0,1,\ldots ,N_t$, with the  step size $\Delta t=T/N_t$. For any function $y(x,t)$, we denote by $y_j^k$ the discrete solution at the node $(x_j,t_k)$, i.e., $y_j^k=y(x_j,t_k)$. We further define the following finite difference and average operators
$$
\delta_ty_j^k = \frac{y_j^{k+1}-y_j^k}{\Delta t}, \qquad A_ty_j^k = \frac{y_j^{k+1}+y_j^k}{2}.
$$

In order to obtain a global energy preserving scheme, the first order spatial derivative $\partial_x$ is  discretized with  centered finite differences, resulting  a skew-symmetric matrix  $D\in\mathbb{R}^{N\times N}$
 under periodic boundary conditions \cite{Gong14}
\begin{equation}\label{matrixD}
D=\frac{1}{2\Delta x}D_N \; , \qquad D_N=
\begin{pmatrix}
0  &  1 &   &   &   & -1\\
-1 &  0 & 1 &   &  &   \\
   & -1 & 0 & 1 &  &  \\
   &    &   & \ddots  &  & \\
   &   &   & -1 & 0 & 1\\
1& &  & & -1 & 0
\end{pmatrix}.
\end{equation}
Then, defining the time-dependent state vector $\bar{\bm{z}}=(z_1;\ldots ;z_N):[0,T]\mapsto\mathbb{R}^{mN}$ with
$$
z_j:=z(x_j,t)=(z_{(1),j},\ldots ,z_{(m),j})^T:[0,T]\mapsto\mathbb{R}^{m}, \qquad j=1,\ldots ,N,
$$
where $z_{(l),j}=z_{(l)}(x_j,t)$, $l=1,\ldots ,m$,
we obtain the following semi-discrete formulation
\begin{equation}\label{mssemi}
K\partial_tz_j + L(D\bar{\bm{z}})_j = \nabla _z S(z_j), \qquad j=1,\ldots, N,
\end{equation}
of the multi-symplectic PDE \eqref{ms}, where we define that $(D\bar{\bm{z}})_j=\sum_{i=1}^ND_{ji}z_i\in\mathbb{R}^{m}$.

\begin{remark}
Note that we distinguish the state components and discrete spatial points by the use of a subscript with or without  parentheses. A subscript $(l)$ with parentheses refers to the $l$th state component, $l=1,\ldots ,m$, whereas a subscript $j$ without parentheses is related to the solution at the discrete spatial point $x_j$, $j=1,\ldots ,N+1$.
\end{remark}

Next, we define the discrete state vector $\bar{\bm{z}}^k=(z_1^k;\ldots ;z_N^k)\in\mathbb{R}^{mN}$, and we apply the global energy preserving AVF method to the semi-discrete system \eqref{mssemi}, to obtain the full discrete system
\begin{align}\label{msfull}
K\delta_tz_j^k + L(DA_t\bar{\bm{z}}^k)_j &= \int_0^1 \nabla _z S((1-\xi )z_j^k+\xi z_j^{k+1})d\xi, & j=1,\ldots, N.
\end{align}

Then,  the discrete  global energy reads as
\begin{equation}\label{denergy}
\varepsilon_h^k=\varepsilon_h(\bar{\bm{z}}^k)=\Delta x\sum_{j=1}^N E_j^k, \qquad E_j^k = S(z_j^k) + (D\bar{\bm{z}}^k)_j^TL_+z_j^k,
\end{equation}
which satisfies the discrete global energy conservation law \cite[Theorem 3.2]{Gong14}, i.e., $\varepsilon_h^{k+1}=\varepsilon_h^k$.

For the  construction of  the ROM with the POD basis functions in Section~\ref{sec:rom}, one needs snapshots of the full-order solutions.
At this point, in order to obtain the FOM, we need to write the spatially component-wise (local) equations \eqref{mssemi} in a compact form with a full set of solutions, i.e., the solution vector in the compact form has to include the solutions at all the spatial nodes $x_j$ in it. On the other hand, the formulation so far relies on the global solution vector $\bar{\bm{z}}=(\bar{\bm{z}}_1;\ldots ;\bar{\bm{z}}_N)$ with the ordering
\begin{equation}\label{ordering1}
\bar{\bm{z}}=(  \underbrace{z_{(1),1},\ldots ,z_{(m),1}}_{\bar{\bm{z}}_1^T}, \ldots , \underbrace{z_{(1),j},\ldots ,z_{(m),j}}_{\bar{\bm{z}}_j^T}, \ldots, \underbrace{z_{(m),N},\ldots ,z_{(m),N}}_{\bar{\bm{z}}_N^T})^T:[0,T]\mapsto\mathbb{R}^{mN},
\end{equation}
where $\bar{\bm{z}}_j=(z_{(1),j},\ldots ,z_{(m),j})^T$ with $z_{(l),j}=z_{(l)}(x_j,t)$. However, for the computational purpose, in the construction of the ROM, a different ordering of the global solution vector is needed.
Let us set the global solution vector $\bm{z}$ in the following form
\begin{equation}\label{ordering}
\bm{z}=(  \underbrace{z_{(1),1},\ldots ,z_{(1),N}}_{\bm{z}_{(1)}^T}, \ldots , \underbrace{z_{(l),1},\ldots ,z_{(l),N}}_{\bm{z}_{(l)}^T}, \ldots, \underbrace{z_{(m),1},\ldots ,z_{(m),N}}_{\bm{z}_{(m)}^T})^T:[0,T]\mapsto\mathbb{R}^{mN},
\end{equation}
where $\bm{z}_{(l)}=(z_{(l),1},\ldots ,z_{(l),N})^T:[0,T]\mapsto\mathbb{R}^{N}$ is the semi-discrete solution vector to the $l$th state component at the discrete spatial nodes.
Let us define the following $mN$-dimensional global system matrices, which are compatible with the ordering \eqref{ordering} of the global solution vector $\bm{z}=(\bm{z}_{(1)};\ldots ;\bm{z}_{(m)})$
\begin{equation}\label{mattens}
\bm{K}=K\otimes I_N\; , \qquad \bm{L}=L\otimes I_N\; , \qquad \bm{D}=I_m\otimes D,
\end{equation}
where $I_N$ and $I_m$ are the $N$-dimensional and $m$-dimensional identity matrices, respectively. We note that by the properties of Kronecker product $\otimes$, the matrices $\bm{K}$, $\bm{L}$ and $\bm{D}$ are all still skew-symmetric matrices.
Finally, we obtain  the FOM as the following ODE system
\begin{equation}\label{dfom}
\bm{K}\dot{\bm{z}} + \bm{L}\bm{D}\bm{z} = \nabla_{\bm{z}} \bm{S}(\bm{z}).
\end{equation}

The full discrete system in matrix-vector form equivalent to \eqref{msfull} is obtained by defining the discrete state vector $\bm{z}^k=(z_{(1)}^k;\ldots ;z_{(m)}^k)\in\mathbb{R}^{mN}$, and applying the AVF method to the semi-discrete system \eqref{dfom}, which yields
\begin{align}\label{msfull2}
\bm{K}\delta_t\bm{z}^k + \bm{L}\bm{D}A_t\bm{z}^k = \int_0^1 \nabla_{\bm{z}} \bm{S}((1-\xi )\bm{z}^k+\xi \bm{z}^{k+1})d\xi.
\end{align}

The discrete global energy \eqref{denergy} can be written equivalently in the following compact form
\begin{equation}\label{vecenergy}
\varepsilon_h^k=\Delta x\left( \sum_{j=1}^N \bm{S}_j(\bm{z}^k) + (\bm{D}\bm{z}^k)^T\bm{L}_+\bm{z}^k\right),
\end{equation}
where we define
\begin{equation}\label{S}
\bm{S}(\bm{z}^k)=(\bm{S}_1(\bm{z}^k);\ldots ;\bm{S}_N(\bm{z}^k))\in\mathbb{R}^{N}, \qquad \bm{S}_j(\bm{z}^k)=S(\bar{\bm{z}}_j),
\end{equation}
with $\bar{\bm{z}}_j=(z_{(1),j},\ldots ,z_{(m),j})^T$ as given by the ordering \eqref{ordering1}, and the matrix $\bm{L}_+$ also satisfies the splitting $\bm{L}=\bm{L}_+-\bm{L}_+^T$.

Multi-symplectic integrators for higher dimensional Hamiltonian PDEs are constructed using tensors with a  proper ordering  of the nodes.
For instance, for a two-dimensional domain, we can use the lexicographic ordering of the nodes $\bm{x}_j=(x_{j_x},y_{j_y})$ on the uniform partition of a spatial domain $[x_L,x_R]\times [y_L,y_R]$ with $x_{j_x}=x_L+(j_x-1)\Delta x$, $y_{j_y}=y_L+(j_y-1)\Delta y$, $j_x=1,\ldots ,N_x+1$, $j_y=1,\ldots ,N_y+1$, with the mesh sizes $\Delta x=(x_R-x_L)/N_x$ and $\Delta y=(y_R-y_L)/N_y$, and with the global degree of freedom $N=N_xN_y$.
Then, the semi-discrete system of the two-dimensional variant of the Hamiltonian PDE \eqref{ms1} reads as
$$
\bm{K}\dot{\bm{z}} + \bm{L}\bm{D}_x\bm{z} + \bm{L}\bm{D}_y\bm{z}  = \nabla_{\bm{z}} \bm{S}(\bm{z}),
$$
where the skew-symmetric matrices $\bm{K},\bm{L}\in\mathbb{R}^{mN\times mN}$ are as in \eqref{mattens}. Using the lexicographic ordering, the finite difference matrices $\bm{D}_x,\bm{D}_y\in\mathbb{R}^{mN\times mN}$ are given by
\begin{equation}\label{2dmat}
\begin{aligned}
\bm{D}_x &= I_m\otimes D_x , & D_x = \left( I_{N_y}\otimes \frac{1}{2\Delta x}D_{N_x}\right) \in\mathbb{R}^{N\times N} , \\
\bm{D}_y &= I_m\otimes D_y, & D_y = \left( \frac{1}{2\Delta y}D_{N_y} \otimes I_{N_x}\right) \in\mathbb{R}^{N\times N},
\end{aligned}
\end{equation}
where the matrices $D_{N_x}$ and $D_{N_y}$ are $N_x$ and $N_y$ dimensional, respectively, difference matrices defined in \eqref{matrixD}, and $\otimes$ denotes the Kronecker product.  The remaining setting can be easily done using the $mN$-dimensional matrices $\bm{K}$, $\bm{L}$, $\bm{D}_x$ and $\bm{D}_y$.

%%%%%%%%%%%%%%%%%%%%%%%%%%%%%%%%%%%%%%%%%%%%%%%%%%%%%%%%%%%%%%%%%%%%%%%%%%%%%%%%%%%%%%%%%%%%%%%%%%%%%%%%%%%%%%%%%%%%%%%%%%%%%%%%%%%%%%%%%%%%%%%%%%%%%%%%%%%%%%%%%%%%%%%%%%%%%%%%%%%%%%%%
\section{Reduced-order model}
\label{sec:rom}

In this section, the construction of global energy preserving ROMs is described.

%%%%%%%%%%%%%%%%%%%%%%%%%%%%%%%%%%%%%%%%%%%%%%%%%%%%%%%%%%%%%%%%%%%%%%%%%%%%%%%%%%%%%%%%%
\subsection{Global energy preserving POD-ROM (P-ROM)}

Let the matrix
$
\mathcal{Z}_{(l)} = [\bm{z}_{(l)}^1 \; \ldots \; \bm{z}_{(l)}^{N_t}] \in\mathbb{R}^{N\times N_t}
$, $l=1,\ldots ,m$,
denotes the snapshot of the solutions of the $l$th state component from the FOM \eqref{dfom}, where  $\bm{z}_{(l)}^k=\bm{z}_{(l)}(t_k)=(z_{(l),1}^k,\ldots ,z_{(l),N}^k)^T\in\mathbb{R}^{N}$ is the solution vector of the $l$th state component at time $t_k$.
Then, we apply the POD to the snapshots $\mathcal{Z}_{(l)}$ to extract the first $n\ll N$ dominant orthogonal POD modes $\{V_{(l)}^1,\ldots ,V_{(l)}^n\}\subset\mathbb{R}^{N}$ forming the matrix of POD modes $V_{(l)}=[V_{(l)}^1\; \cdots \; V_{(l)}^n]\in\mathbb{R}^{N\times n}$. The POD modes are computed so that the projection error on the snapshot matrix $\mathcal{Z}_{(l)}$ is minimized, i.e., they are given by the solution of the minimization problem
\begin{equation*}
\min_{V_{(l)}\in \mathbb{R}^{ N\times n}}||\mathcal{Z}_{(l)}-V_{(l)}V_{(l)}^T\mathcal{Z}_{(l)} ||_F^2 , \qquad l=1,\ldots ,m,
\end{equation*}
where $\|\cdot\|_F$ denotes the Frobenius norm.
The POD modes are computed from the the singular value decomposition (SVD) of the snapshot matrix $\mathcal{Z}_{(l)}$. The POD modes $\{V_{(l)}^1,\ldots ,V_{(l)}^n\}$ are the left singular vectors of the snapshot matrix $\mathcal{Z}_{(l)}$ corresponding to the first $n$ singular values in descending order.
Once the POD modes are obtained, the reduced-order approximation $\widehat{\bm{z}}_{(l)}:[0,T]\mapsto \mathbb{R}^{N}$ to the full-order solution $\bm{z}_{(l)}:[0,T]\mapsto \mathbb{R}^{N}$ can be obtained from the reduced space by the linear combination of the POD modes
\begin{equation}\label{rapp}
\bm{z}_{(l)}(t) \approx \widehat{\bm{z}}_{(l)}(t) = \sum_{r=1}^n \alpha_{(l),r}(t)V_{(l)}^r=V_{(l)}\bm{\alpha}_{(l)}(t),
\end{equation}
where $\bm{\alpha}_{(l)}=(\alpha_{(l),1},\ldots ,\alpha_{(l),n})^T:[0,T]\mapsto \mathbb{R}^{n}$ is the vector of unknown coefficients (reduced coefficients) to be determined through the reduced system. We note that the approximation \eqref{rapp} holds separately for each state component. In order to construct a reduced-order approximation to the global solution of the FOM, we define the block diagonal global matrix of POD modes $\bm{V}$, given by
$$
\bm{V} =
\begin{pmatrix}
V_{(1)} & & & \\
 &  V_{(2)} &  & \\
&  &  \ddots  &  \\
& & & V_{(m)}
\end{pmatrix}\in\mathbb{R}^{mN\times mn}.
$$
We can write the following relations between the $mN$-dimensional global reduced-order approximation $\widehat{\bm{z}}=(\widehat{\bm{z}}_{(1)};\ldots ;\widehat{\bm{z}}_{(m)})$  and the $mn$-dimensional vector $\bm{\alpha}=(\bm{\alpha}_{(1)};\ldots ;\bm{\alpha}_{(m)})$ of global reduced coefficients
\begin{equation}\label{rapp2}
\bm{z}(t) \approx \widehat{\bm{z}}(t)=\bm{V}\bm{\alpha}(t), \qquad \bm{\alpha}(t) = \bm{V}^T\widehat{\bm{z}}(t),
\end{equation}
where we use the orthogonality property of the global matrix $\bm{V}$ of POD modes, i.e., $\bm{V}^T\bm{V}=I_{mn}$.

Inserting the first approximation in \eqref{rapp2} into the FOM \eqref{dfom}, and projecting  onto the reduced space by multiplying from left by the global matrix $\bm{V}^T$ of POD modes, we get the reduced-order system
\begin{equation}\label{drom}
\bm{V}^T\bm{K}\bm{V}\dot{\bm{\alpha}} + \bm{V}^T\bm{L}\bm{D}\bm{V}\bm{\alpha} = \bm{V}^T\nabla_{\bm{z}} \bm{S}(\bm{V}\bm{\alpha}).
\end{equation}
In order to preserve the skew-symmetric structure in the reduced-order setting, as in the FOM \eqref{dfom}, we simply insert $\bm{V}\bm{V}^T$ between the matrices $\bm{L}$ and $\bm{D}$ in the reduced-order system\eqref{drom}, and we obtain the following reduced-order ODE system
\begin{equation}\label{drom2}
\widehat{\bm{K}}\dot{\bm{\alpha}} + \widehat{\bm{L}}\widehat{\bm{D}}\bm{\alpha} = \nabla_{\bm{\alpha}} \widehat{\bm{S}}(\bm{\alpha}),
\end{equation}
where, $\widehat{\bm{K}}=\bm{V}^T\bm{K}\bm{V}$, $\widehat{\bm{L}}=\bm{V}^T\bm{L}\bm{V}$ and $\widehat{\bm{D}}=\bm{V}^T\bm{D}\bm{V}$ are the $mn$-dimensional reduced skew-symmetric matrices, and the reduced gradient vector is given by
\begin{equation}\label{rg}
\nabla_{\bm{\alpha}} \widehat{\bm{S}}(\bm{\alpha})=\nabla_{\bm{\alpha}} \bm{S}(\bm{V}\bm{\alpha})=\bm{V}^T\nabla_{\bm{z}} \bm{S}(\bm{V}\bm{\alpha}).
\end{equation}

In the sequel, we call the ROM \eqref{drom2} as P-ROM.
Similar to the FOM, the P-ROM \eqref{drom2} is also integrated in time with the global energy preserving AVF method, which yields the full discrete reduced system
\begin{align}\label{romfull}
\widehat{\bm{K}}\delta_t\bm{\alpha}^k + \widehat{\bm{L}}\widehat{\bm{D}}A_t\bm{\alpha}^k = \int_0^1 \nabla_{\bm{\alpha}} \widehat{\bm{S}}((1-\xi )\bm{\alpha}^k+\xi \bm{\alpha}^{k+1})d\xi,
\end{align}
where $\bm{\alpha}^k=\bm{\alpha}(t_k)$ is the vector of reduced coefficients at time $t_k$.

\begin{theorem}\label{thm1}
For the full-rank reduced-order approximation $\widehat{\bm{z}}^k=\widehat{\bm{z}}(t_k)=\bm{V}\bm{\alpha}^k$ obtained with the full discrete P-ROM \eqref{romfull}, and with the choice $L_+=(1/2)L$, the discrete GCL holds, i.e., $\widehat{\varepsilon}_h^{k+1}=\widehat{\varepsilon}_h^k$, with the discrete reduced global energy
\begin{equation}\label{vecenergy2}
\widehat{\varepsilon}_h^k=\widehat{\varepsilon}_h(\bm{\alpha}^k)=\Delta x\left( \sum_{j=1}^N \widehat{\bm{S}}_j(\bm{\alpha}^k) + (\widehat{\bm{D}}\bm{\alpha}^k)^T\widehat{\bm{L}}_+\bm{\alpha}^k\right),
\end{equation}
where $\widehat{\bm{L}}_+=\bm{V}^T\bm{L}_+\bm{V}$ and $\widehat{\bm{S}}_j(\bm{\alpha}^k)=\bm{S}_j(\bm{V}\bm{\alpha}^k)=\bm{S}_j(\widehat{\bm{z}}^k)$.

\begin{proof}
The proof relies on showing that $\delta_t\widehat{\varepsilon}_h^k=0$ which is equivalent to the identity $\widehat{\varepsilon}_h^{k+1}=\widehat{\varepsilon}_h^k$.
We take the inner product with $\delta_t\bm{\alpha}^k$ on both sides of the full discrete P-ROM \eqref{romfull} to have the system
\begin{align}\label{prf1}
\left(\delta_t\bm{\alpha}^k\right)^T\widehat{\bm{L}}\widehat{\bm{D}}A_t\bm{\alpha}^k = \left(\delta_t\bm{\alpha}^k\right)^T\int_0^1 \nabla_{\bm{\alpha}} \widehat{\bm{S}}((1-\xi )\bm{\alpha}^k+\xi \bm{\alpha}^{k+1})d\xi,
\end{align}
where we use the fact that $\left(\delta_t\bm{\alpha}^k\right)^T\widehat{\bm{K}}\delta_t\bm{\alpha}^k=0$ because of the skew-symmetric matrix $\widehat{\bm{K}}$. Using the properties of the integral and the definition \eqref{rg} of the reduced gradient, we can write the right hand side integral term in \eqref{prf1} as
\begin{align*}
\left(\delta_t\bm{\alpha}^k\right)^T\int_0^1 &\nabla_{\bm{\alpha}} \widehat{\bm{S}}((1-\xi )\bm{\alpha}^k+ \xi \bm{\alpha}^{k+1})d\xi \\
&= \frac{1}{\Delta t}(\bm{\alpha}^{k+1}-\bm{\alpha}^k)^T\int_0^1 \bm{V}^T\nabla_{\bm{z}} \bm{S}((1-\xi )\bm{V}\bm{\alpha}^k+\xi \bm{V}\bm{\alpha}^{k+1})d\xi\\
&= \frac{1}{\Delta t}(\bm{V}\bm{\alpha}^{k+1}-\bm{V}\bm{\alpha}^k)^T\int_0^1 \nabla_{\bm{z}} \bm{S}((1-\xi )\bm{V}\bm{\alpha}^k+\xi \bm{V}\bm{\alpha}^{k+1})d\xi\\
&= \frac{1}{\Delta t}(\widehat{\bm{z}}^{k+1}-\widehat{\bm{z}}^k)^T\int_0^1 \nabla_{\bm{z}} \bm{S}((1-\xi )\widehat{\bm{z}}^k+\xi \widehat{\bm{z}}^{k+1})d\xi.
\end{align*}
For the solution vector $\widehat{\bar{\bm{z}}}^k$ which is the reordered form of the vector $\widehat{\bm{z}}^k$ with the ordering \eqref{ordering1},  the above identity yields
\begin{align*}
\left(\delta_t\bm{\alpha}^k\right)^T\int_0^1 \nabla_{\bm{\alpha}} \widehat{\bm{S}}((1-\xi )\bm{\alpha}^k &+ \xi \bm{\alpha}^{k+1})d\xi \\
&= \frac{1}{\Delta t}(\widehat{\bar{\bm{z}}}^{k+1}-\widehat{\bar{\bm{z}}}^k)^T\int_0^1 \nabla_{\bar{\bm{z}}} \bm{S}((1-\xi )\widehat{\bar{\bm{z}}}^k+\xi \widehat{\bar{\bm{z}}}^{k+1})d\xi\\
&= \frac{1}{\Delta t}\sum_{j=1}^N \int_0^1 \frac{d}{d\xi} \bm{S}_j((1-\xi )\widehat{\bm{z}}^k+\xi \widehat{\bm{z}}^{k+1})d\xi\\
&= \delta_t\sum_{j=1}^N\bm{S}_j(\widehat{\bm{z}}^k)\\
&= \delta_t\sum_{j=1}^N \widehat{\bm{S}}_j(\bm{\alpha}^k).
\end{align*}
Using the above identity and utilizing the splitting $\widehat{\bm{L}}=\widehat{\bm{L}}_+-\widehat{\bm{L}}_+^T$ of the skew-symmetric matrix $\widehat{\bm{L}}$, we get from \eqref{prf1}
\begin{equation}\label{deltaS}
\delta_t\sum_{j=1}^N \widehat{\bm{S}}_j(\bm{\alpha}^k) = \left(\delta_t\bm{\alpha}^k\right)^T\widehat{\bm{L}}_+\widehat{\bm{D}}A_t\bm{\alpha}^k - \left(\widehat{\bm{D}}A_t\bm{\alpha}^k\right)^T\widehat{\bm{L}}_+\delta_t\bm{\alpha}^k.
\end{equation}
According to commutative law and discrete Leibniz rule
$$
\delta_t(f\cdot g)^k = \delta_tf^k\cdot A_tg^k + A_tf^k\cdot \delta_tg^k,
$$
we can write the identity
\begin{equation}\label{leib}
\delta_t((\widehat{\bm{D}}\bm{\alpha}^k)^T\widehat{\bm{L}}_+\bm{\alpha}^k) = (\widehat{\bm{D}}\delta_t\bm{\alpha}^k)^T\widehat{\bm{L}}_+A_t\bm{\alpha}^k + (\widehat{\bm{D}}A_t\bm{\alpha}^k)^T\widehat{\bm{L}}_+\delta_t\bm{\alpha}^k.
\end{equation}

Using the identity \eqref{leib} and the equation \eqref{deltaS}, we get for the discrete reduced global energy \eqref{vecenergy2}
\begin{equation}\label{prfen2}
\begin{aligned}
\delta_t\widehat{\varepsilon}_h^k &=\Delta x\left( \delta_t\sum_{j=1}^N \widehat{\bm{S}}_j(\bm{\alpha}^k) + (\widehat{\bm{D}}\delta_t\bm{\alpha}^k)^T\widehat{\bm{L}}_+A_t\bm{\alpha}^k + (\widehat{\bm{D}}A_t\bm{\alpha}^k)^T\widehat{\bm{L}}_+\delta_t\bm{\alpha}^k\right)\\
&= \Delta x\left( \left(\delta_t\bm{\alpha}^k\right)^T\widehat{\bm{L}}_+\widehat{\bm{D}}A_t\bm{\alpha}^k + (\widehat{\bm{D}}\delta_t\bm{\alpha}^k)^T\widehat{\bm{L}}_+A_t\bm{\alpha}^k \right)\\
&= \Delta x \left(\delta_t\bm{\alpha}^k\right)^T\left( \widehat{\bm{L}}_+\widehat{\bm{D}} +\widehat{\bm{D}}^T\widehat{\bm{L}}_+ \right) A_t\bm{\alpha}^k\\
&= \Delta x \left(\frac{\bm{\alpha}^{k+1}-\bm{\alpha}^{k}}{\Delta t}\right)^T\left( \widehat{\bm{L}}_+\widehat{\bm{D}} +\widehat{\bm{D}}^T\widehat{\bm{L}}_+ \right) \left(\frac{\bm{\alpha}^{k+1}+\bm{\alpha}^{k}}{2}\right)\\
&= \frac{\Delta x}{2\Delta t}\left( \left(\bm{\alpha}^{k+1}\right)^T\left( \widehat{\bm{L}}_+\widehat{\bm{D}} +\widehat{\bm{D}}^T\widehat{\bm{L}}_+ \right) \bm{\alpha}^{k+1}  -  \left(\bm{\alpha}^{k}\right)^T\left( \widehat{\bm{L}}_+\widehat{\bm{D}} +\widehat{\bm{D}}^T\widehat{\bm{L}}_+ \right) \bm{\alpha}^{k} \right).
\end{aligned}
\end{equation}
Then, the discrete GCL $\delta_t\widehat{\varepsilon}_h^k=0$ can be obtained if the matrix $\left( \widehat{\bm{L}}_+\widehat{\bm{D}} +\widehat{\bm{D}}^T\widehat{\bm{L}}_+ \right)$ is skew-symmetric. Indeed, for the choice $\widehat{\bm{L}}_+=(1/2)\widehat{\bm{L}}$, it is true by the skew-symmetry of the matrix $\widehat{\bm{L}}$, since we have then
\begin{align*}
\left( \widehat{\bm{L}}_+\widehat{\bm{D}} +\widehat{\bm{D}}^T\widehat{\bm{L}}_+ \right)^T &= \frac{1}{2} \left( \widehat{\bm{L}}\widehat{\bm{D}} +\widehat{\bm{D}}^T\widehat{\bm{L}} \right)^T \\
&= \frac{1}{2} \left( \widehat{\bm{D}}^T\widehat{\bm{L}}^T + \widehat{\bm{L}}^T\widehat{\bm{D}}  \right)  \\
&= \frac{1}{2} \left( -\widehat{\bm{D}}^T\widehat{\bm{L}} - \widehat{\bm{L}}\widehat{\bm{D}} \right)\\
&= - \frac{1}{2} \left( \widehat{\bm{L}}\widehat{\bm{D}} +\widehat{\bm{D}}^T\widehat{\bm{L}} \right)\\
&= - \left( \widehat{\bm{L}}_+\widehat{\bm{D}} +\widehat{\bm{D}}^T\widehat{\bm{L}}_+ \right).
\end{align*}

\end{proof}
\end{theorem}

%%%%%%%%%%%%%%%%%%%%%%%%%%%%%%%%%%%%%%%%%%%%%%%%%%%%%%%%%%%%%%%%%%%%%%%%%%%%%%%%%%%%%%%%%
\subsection{Global energy preserving POD-DEIM-ROM (PD-ROM)}

The computation of the linear left-hand side terms in the P-ROM \eqref{drom2} scales with the reduced dimension $mn$, while the nonlinear right-hand side term scales with the full dimension $mN$.
Applying hyper-reduction techniques such as EIM and DEIM \cite{chaturantabut10nmr,Gugercin16a,Barrault04}, the offline-online computations are decomposed, i.e., the ROMs are independent of the FOM dimension.
The full dimensional gradient vector $\nabla_{\bm{z}} \bm{S}(\bm{V}\bm{\alpha})$ is approximated with the DEIM  through the interpolation onto an empirical basis
\begin{equation*}
\nabla_{\bm{z}} \bm{S}(\bm{V}\bm{\alpha}) \approx \bm{\Phi} \mathsf{c}(t),
\end{equation*}
where $\bm{\Phi} = [\bm{\Phi}_{1}  \cdots  \bm{\Phi}_{\tilde{n}}] \in \mathbb{R}^{mN \times \tilde{n}}$ is a low dimensional ($\tilde{n}\ll N$) basis matrix and $\mathsf{c}(t) : [0,T] \mapsto \mathbb{R}^{\tilde{n}}$ is the vector of time-dependent coefficients to be determined.
The coefficient vector $\mathsf{c}(t)$ can be calculated from a projected system $\mathsf{P}^T\nabla_{\bm{z}} \bm{S}(\bm{V}\bm{\alpha})=\mathsf{P}^T\bm{\Phi} \mathsf{c}(t)$ for a selection (permutation) matrix $\mathsf{P}\in \mathbb{R}^{mN \times \tilde{n}}$ so that the matrix $\mathsf{P}^T\bm{\Phi}$ is non-singular.
Once the selection matrix $\mathsf{P}$ is determined, the coefficient vector $\mathsf{c}(t)$ is calculated, which yields the approximation
\begin{equation}\label{deimapp}
\nabla_{\bm{z}} \bm{S}(\bm{V}\bm{\alpha}) \approx \bm{W} \nabla_{\bm{z}} \bm{S}_D(\bm{V}\bm{\alpha}),
\end{equation}
where the constant matrix $\bm{W}=\bm{\Phi}\left(\mathsf{P}^T\bm{\Phi} \right)^{-1}\in\mathbb{R}^{mN \times \tilde{n}}$ can be precomputed in the offline stage, and the DEIM approximated nonlinear term is given as
\begin{equation*}
\nabla_{\bm{z}}\bm{S}_D(\bm{V}\bm{\alpha})=\mathsf{P}^T\nabla_{\bm{z}}\bm{S}(\bm{V}\bm{\alpha}):[0,T]\mapsto \mathbb{R}^{\tilde{n}},
\end{equation*}
which does not scale with the full dimension.
The DEIM approximation is constructed  by selecting $\tilde{n}$ entries of the nonlinear vector $\nabla_{\bm{z}}\bm{S}(\bm{V}\bm{\alpha})$ among $mN$ entries.
The selection matrix  $\mathsf{P}$ is computed by the QDEIM algorithm \cite{gugercin16}, which has better stability properties and more accurate than the DEIM \cite{chaturantabut10nmr}.
The DEIM basis matrix $\bm{\Phi}$ is obtained with the SVD  of the nonlinear snapshot matrix $\mathcal{S} = [ \nabla_{\bm{z}}\bm{S}(\bm{V}\bm{\alpha}^1),\cdots , \nabla_{\bm{z}}\bm{S}(\bm{V}\bm{\alpha}^{N_t})]\in\mathbb{R}^{mN\times N_t}$.

Inserting the DEIM approximation \eqref{deimapp}, the P-ROM \eqref{drom2} takes the form
\begin{equation}\label{drom3}
\widehat{\bm{K}}\dot{\bm{\alpha}} + \widehat{\bm{L}}\widehat{\bm{D}}\bm{\alpha} = \nabla_{\bm{\alpha}} \widehat{\bm{S}}_D(\bm{\alpha}),
\end{equation}
with the  reduced nonlinear term
\begin{equation}\label{deimapp2}
\nabla_{\bm{\alpha}} \widehat{\bm{S}}_D(\bm{\alpha})=\bm{V}^T\bm{W}\nabla_{\bm{z}} \bm{S}_D(\bm{V}\bm{\alpha}), \qquad \nabla_{\bm{z}} \bm{S}_D(\bm{V}\bm{\alpha})=\mathsf{P}^T\nabla_{\bm{z}} \bm{S}(\bm{V}\bm{\alpha}).
\end{equation}
In the sequel, we call the ROM \eqref{drom3} as POD-DEIM-ROM (PD-ROM) in order to distinguish the ROMs with and without DEIM approximation.
The PD-ROM \eqref{drom3} is again integrated in time with the global energy preserving AVF method
\begin{align}\label{romfullD}
\widehat{\bm{K}}\delta_t\bm{\alpha}^k + \widehat{\bm{L}}\widehat{\bm{D}}A_t\bm{\alpha}^k = \int_0^1 \nabla_{\bm{\alpha}} \widehat{\bm{S}}_D((1-\xi )\bm{\alpha}^k+\xi \bm{\alpha}^{k+1})d\xi.
\end{align}

The PD-ROM \eqref{drom3} has the same form as the P-ROM \eqref{drom2}. However, the discrete reduced global energy \eqref{vecenergy2} is conserved approximately.

\begin{theorem}
For the reduced approximation $\widehat{\bm{z}}^k=\widehat{\bm{z}}(t_k)=\bm{V}\bm{\alpha}^k$ obtained through the full discrete PD-ROM \eqref{romfullD} with DEIM approximation, and with the choice $L_+=(1/2)L$, the discrete reduced global energy $\widehat{\varepsilon}_h^k$ in \eqref{vecenergy2} is approximately conserved. The error bound for the discrete reduced global energy is given by
\begin{equation}\label{errbound}
\|\delta_t\widehat{\varepsilon}_h^k\| \le \sum_{j=1}^N  \Delta x\|(\mathsf{P}^T\bm{\Phi})^{-1}\|\|(\bm{I}-\bm{\Phi}\bm{\Phi}^T)\| \|\delta_t \widehat{\bm{S}}_j(\bm{\alpha}^k)\|.
\end{equation}

\begin{proof}
We take the inner product with $\delta_t\bm{\alpha}^k$ on both sides of the full discrete PD-ROM \eqref{romfullD} to have the system
\begin{align}\label{prf2}
\left(\delta_t\bm{\alpha}^k\right)^T\widehat{\bm{L}}\widehat{\bm{D}}A_t\bm{\alpha}^k = \left(\delta_t\bm{\alpha}^k\right)^T\int_0^1 \nabla_{\bm{\alpha}} \widehat{\bm{S}}_D((1-\xi )\bm{\alpha}^k+\xi \bm{\alpha}^{k+1})d\xi,
\end{align}
where we again use the fact that $\left(\delta_t\bm{\alpha}^k\right)^T\widehat{\bm{K}}\delta_t\bm{\alpha}^k=0$ because of the skew-symmetric matrix $\widehat{\bm{K}}$.
Using the properties of the integral and the definition \eqref{deimapp2} of the DEIM approximated nonlinear term, we can write the right hand side integral term in \eqref{prf2} as
\begin{align*}
\left(\delta_t\bm{\alpha}^k\right)^T\int_0^1 &\nabla_{\bm{\alpha}} \widehat{\bm{S}}_D((1-\xi )\bm{\alpha}^k+\xi \bm{\alpha}^{k+1})d\xi \\
&= \frac{1}{\Delta t}(\bm{\alpha}^{k+1}-\bm{\alpha}^k)^T\int_0^1 \bm{V}^T\bm{W}\nabla_{\bm{z}} \bm{S}_D((1-\xi )\bm{V}\bm{\alpha}^k+\xi \bm{V}\bm{\alpha}^{k+1})d\xi\\
&= \frac{1}{\Delta t}(\bm{V}\bm{\alpha}^{k+1}-\bm{V}\bm{\alpha}^k)^T\int_0^1 \bm{W}\nabla_{\bm{z}} \bm{S}_D((1-\xi )\bm{V}\bm{\alpha}^k+\xi \bm{V}\bm{\alpha}^{k+1})d\xi\\
&= \frac{1}{\Delta t}(\widehat{\bm{z}}^{k+1}-\widehat{\bm{z}}^k)^T\int_0^1 \bm{W}\nabla_{\bm{z}} \bm{S}_D((1-\xi )\widehat{\bm{z}}^k+\xi \widehat{\bm{z}}^{k+1})d\xi.
\end{align*}
Adding and subtracting the term $\nabla_{\bm{z}} \bm{S}((1-\xi )\widehat{\bm{z}}^k+\xi \widehat{\bm{z}}^{k+1})$ in the integrand of the above identity, and using the result of the same step in the proof of Theorem~\ref{thm1}, we obtain the equation
\begin{align*}
&\left(\delta_t\bm{\alpha}^k\right)^T  \int_0^1 \nabla_{\bm{\alpha}} \widehat{\bm{S}}_D((1-\xi )\bm{\alpha}^k+\xi \bm{\alpha}^{k+1})d\xi = \delta_t\sum_{j=1}^N \widehat{\bm{S}}_j(\bm{\alpha}^k) \\
& - \frac{1}{\Delta t}(\widehat{\bm{z}}^{k+1}-\widehat{\bm{z}}^k)^T\int_0^1 \left( \nabla_{\bm{z}} \bm{S}((1-\xi )\widehat{\bm{z}}^k+\xi \widehat{\bm{z}}^{k+1}) - \bm{W}\nabla_{\bm{z}} \bm{S}_D((1-\xi )\widehat{\bm{z}}^k+\xi \widehat{\bm{z}}^{k+1}) \right) d\xi.
\end{align*}
Following the similar steps in the proof of Theorem~\ref{thm1}, the reduced global energy is obtained in the following form
\begin{equation}\label{prfen3}
\delta_t\widehat{\varepsilon}_h^k = - \frac{\Delta x}{\Delta t}(\widehat{\bm{z}}^{k+1}-\widehat{\bm{z}}^k)^T\int_0^1 \left( \nabla_{\bm{z}} \bm{S}((1-\xi )\widehat{\bm{z}}^k+\xi \widehat{\bm{z}}^{k+1}) - \bm{W}\nabla_{\bm{z}} \bm{S}_D((1-\xi )\widehat{\bm{z}}^k+\xi \widehat{\bm{z}}^{k+1}) \right) d\xi.
\end{equation}
Using the DEIM approximation error \cite[Lemma 3.2]{chaturantabut10nmr}, the integrand  in \eqref{prfen3} is bounded by
\begin{equation}\label{deimerror}
\begin{aligned}
\|\nabla_{\bm{z}} \bm{S}((1-\xi )\widehat{\bm{z}}^k+ &\xi \widehat{\bm{z}}^{k+1}) - \bm{W}\nabla_{\bm{z}} \bm{S}_D((1-\xi )\widehat{\bm{z}}^k+\xi \widehat{\bm{z}}^{k+1})\| \\
&\le \|(\mathsf{P}^T\bm{\Phi})^{-1}\|\|(\bm{I}-\bm{\Phi}\bm{\Phi}^T)\|\|\nabla_{\bm{z}} \bm{S}((1-\xi )\widehat{\bm{z}}^k+\xi \widehat{\bm{z}}^{k+1})\|,
\end{aligned}
\end{equation}
where $\|\cdot\|$ denotes the Euclidean $2$-norm. Finally, using the error bound \eqref{deimerror} inside \eqref{prfen3}, we obtain the error bound for the DEIM approximated reduced global energy
\begin{align*}
\|\delta_t\widehat{\varepsilon}_h^k\| &\le \Delta x \|(\mathsf{P}^T\bm{\Phi})^{-1}\|\|(\bm{I}-\bm{\Phi}\bm{\Phi}^T)\| \left\| \frac{1}{\Delta t}(\widehat{\bm{z}}^{k+1}-\widehat{\bm{z}}^k)^T\int_0^1 \nabla_{\bm{z}} \bm{S}((1-\xi )\widehat{\bm{z}}^k+\xi \widehat{\bm{z}}^{k+1})  d\xi\right\|\\
&\le  \sum_{j=1}^N  \Delta x\|(\mathsf{P}^T\bm{\Phi})^{-1}\|\|(\bm{I}-\bm{\Phi}\bm{\Phi}^T)\| \|\delta_t \widehat{\bm{S}}_j(\bm{\alpha}^k)\|.
\end{align*}

\end{proof}
\end{theorem}

\begin{remark}
We note that \eqref{errbound} provides a well defined estimate for the DEIM approximated reduced global energy since the columns of the matrix $\bm{\Phi}$ are orthonormal so that the term $\|(\mathsf{P}^T\bm{\Phi})^{-1}\|$ is of moderate size. Moreover, the matrix $\bm{\Phi}\bm{\Phi}^T$ converges to the identity matrix by increasing number of DEIM basis modes $\tilde{n}$.
\end{remark}

%%%%%%%%%%%%%%%%%%%%%%%%%%%%%%%%%%%%%%%%%%%%%%%%%%%%%%%%%%%%%%%%%%%%%%%%%%%%%%%%%%%%%%%%%%%%
%%%%%%%%%%%%%%%%%%%%%%%%%%%%%%%%%%%%%%%%%%%%%%%%%%%%%%%%%%%%%%%%%%%%%%%%%%%%%%%%%%%%%%%%%%%%
\section{Numerical results}
\label{sec:num}

In this section, we demonstrate the accuracy and computational efficiency of the ROMs for one- and two-dimensional multi-symplectic Hamiltonian PDEs \eqref{ms1}: the one-dimensional KdV and NLS equations, and the two-dimensional NLS and ZK equations, where the ZK equation is a two-dimensional generalization of the KdV equation. All the models are prescribed with periodic boundary conditions.
For the one-dimensional KdV and NLS equations, the construction of the FOMs and  ROMs are given together with the discrete global energies in a detailed form.
For the two-dimensional NLS and ZK equations, we present only the discrete global energies and the solution of the FOMs and ROMs.

In all the numerical tests, we consider solitons with known exact solutions.
The accuracy of the FOM solutions and reduced approximations are measured using  the relative error
\begin{equation}\label{solerr}
\text{E}_{\text{sol}} = \frac{\|Z_h-Z_E\|_F}{\|Z_E\|_F},
\end{equation}
in the Frobenius norm, where $Z_h\in\mathbb{R}^{N\times N_t}$ is the matrix of numerical solutions (by FOM or ROM) with $(Z_h)_{j,k}=\bm{z}_j^k$, and $Z_E\in\mathbb{R}^{N\times N_t}$ is the matrix of exact solutions with $(Z_E)_{j,k}=z(x_j,t_k)$.

The relative shape error is defined by
\begin{equation}\label{shape}
\text{E}_{\text{shape}}:=\min_k \frac{\| \bm{z}^{N_t}-z(\cdot,t_k)\|_2^2}{\| z(\cdot,T)\|_2^2} ,
\end{equation}
in the $2$-norm, where $\bm{z}^{N_t}$ is the numerical solution at the final time $t_{N_t}=T$, and $z(\cdot,t_k)$ is the exact solution at time $t_k$.
Preservation of the discrete global energy is measured with relative energy error
\begin{equation}\label{enerr}
\text{E}_{\text{energy}} = \max_k \frac{|\varepsilon_h^k-\varepsilon_h^0|}{|\varepsilon_h^0|}.
\end{equation}

The number of POD and DEIM modes can be determined using different criteria.
Here, we determine both numbers by selecting the smallest index $i$ such that the normalized singular values of the snapshot matrix satisfy $\sigma_i/\sigma_1 < \tau$ for a given tolerance $\tau >0$.
In our simulations, we take $\tau \approx 10^{-3}$ to determine the number of POD modes.
The nonlinearity should be approximated more accurately, therefore the DEIM modes are determined with $\tau \approx 10^{-5}$.
In the case of the one-dimensional KdV equation, ROMs are computed for an increasing number of POD modes to illustrate the solution and shape errors, the behavior of the solitons, and preservation of the discrete global energies.

All the simulations are performed on a machine with Intel CoreTM i7 2.5 GHz 64 bit CPU, 16 GB RAM, Windows 10, using 64 bit MatLab R2014.

%%%%%%%%%%%%%%%%%%%%%%%%%%%%%%%%%%%%%%%%%%%%%%%%%%%%%%%%%%%%%%%%%%%%%%%%%%%%%%%%%%%%%%%%%%%%
\subsection{Korteweg-de Vries equation}
\label{ex:kdv}

The full discrete FOM \eqref{msfull2} of the multi-symplectic KdV equation \eqref{kdvms1} is given by
\begin{equation}\label{fullkdvms}
\begin{aligned}
\frac{1}{2}\delta_t\bm{u}^k+DA_t\bm{w}^k&=0,\\
-\frac{1}{2}\delta_t\bm{\phi}^k-\gamma DA_t\bm{v}^k &= -A_t\bm{w}^k + \int_0^1 \bm{f}((1-\xi)\bm{u}^k + \xi\bm{u}^{k+1})d\xi ,\\
\gamma DA_t\bm{u}^k&=A_t\bm{v}^k,\\
-DA_t\bm{\phi}^k &= -A_t\bm{u}^k,
\end{aligned}
\end{equation}
with $\bm{\phi}=\bm{z}_{(1)}$, $\bm{u}=\bm{z}_{(2)}$, $\bm{v}=\bm{z}_{(3)}$, $\bm{w}=\bm{z}_{(4)}$, and $\bm{f}(\bm{u})=\eta \bm{u}^2/2$
is the nonlinear term.
Eliminating the auxiliary variables, the system \eqref{fullkdvms} is solved for the unknown vector $\bm{u}$ only
\begin{equation}\label{fullkdvms2}
\delta_t\bm{u}^k+\gamma^2 D^3A_t\bm{u}^k = -D\int_0^1 \bm{f}((1-\xi)\bm{u}^k + \xi\bm{u}^{k+1})d\xi .
\end{equation}
The discrete global energy \eqref{vecenergy} corresponding to the global continuous energy \eqref{kdvenergy1} is given by
$$
\varepsilon_h^k = \Delta x \sum_{j=1}^N \left( \frac{\eta}{6}(\bm{u}_j^k)^3 - \frac{\gamma^2}{2}(D\bm{u}^k)_j^2\right).
$$

The same POD basis matrix $V\in\mathbb{R}^{N\times n}$ in \eqref{rapp} is taken for each state
$$
\bm{\phi}(t) \approx V\bm{\alpha}_{(1)}(t) , \quad \bm{u}(t)\approx V\bm{\alpha}_{(2)}(t), \quad \bm{v}(t)\approx V\bm{\alpha}_{(3)}(t), \quad \bm{w}(t)\approx V\bm{\alpha}_{(4)}(t),
$$
where the vectors $\bm{\alpha}_{(i)}(t): [0,T]\mapsto \mathbb{R}^{n}$ are the reduced coefficients, $i=1,\ldots ,4$.
By this choice, the reduced solutions can be computed  by solving the reduced-order system for a single unknown vector as for the FOM \eqref{fullkdvms2}.

The POD basis matrix $V$ is computed with  the snapshot matrix
$$
\mathcal{U} = [\bm{u}^1, \; \ldots \;,  \bm{u}^{N_t} \; \gamma D\bm{u}^1,\; \ldots \;, \gamma D\bm{u}^{N_t} ] \in\mathbb{R}^{N\times 2N_t},
$$
by concatenating the states $u$ and $v=\gamma u_x$.
Using  $V^TV=I_n$, the  discrete P-ROM \eqref{romfull} for the KdV equation read as
\begin{equation}\label{romfullkdv}
\begin{aligned}
\frac{1}{2}\delta_t\bm{\alpha}_{(2)}^k+\widehat{D}A_t\bm{\alpha}_{(4)}^k&=0,\\
-\frac{1}{2}\delta_t\bm{\alpha}_{(1)}^k-\gamma \widehat{D}A_t\bm{\alpha}_{(3)}^k &= -A_t\bm{\alpha}_{(4)}^k + V^T\int_0^1 \bm{f}((1-\xi)V\bm{\alpha}_{(2)}^k + \xi V\bm{\alpha}_{(2)}^{k+1})d\xi ,\\
\gamma \widehat{D}A_t\bm{\alpha}_{(2)}^k&=A_t\bm{\alpha}_{(3)}^k,\\
-\widehat{D}A_t\bm{\alpha}_{(1)}^k &= -A_t\bm{\alpha}_{(2)}^k,
\end{aligned}
\end{equation}
where $\widehat{D}=V^TDV\in\mathbb{R}^{n\times n}$ denotes the reduced skew-symmetric matrix.
Similar to the full discrete FOM, we eliminate the coefficient vectors $\bm{\alpha}_{(i)}$ of the auxiliary variables, $i=1,3,4$, so that the system \eqref{romfullkdv} can be solved for  the unknown vector $\bm{\alpha}_{(2)}$ only
\begin{equation}\label{romfullkdv2}
\delta_t\bm{\alpha}_{(2)}^k+\gamma^2 \widehat{D}^3A_t\bm{\alpha}_{(2)}^k = -\widehat{D}V^T\int_0^1 \bm{f}((1-\xi)V\bm{\alpha}_{(2)}^k + \xi V\bm{\alpha}_{(2)}^{k+1})d\xi .
\end{equation}

The PD-ROM is obtained by approximating the nonlinear vector $\bm{f}(\bm{u})=\eta \bm{u}^2/2$ from the column space of the nonlinear snapshot matrix
$
\mathcal{S} = [\bm{f}^1 \; \ldots \; \bm{f}^{N_t}  ] \in\mathbb{R}^{N\times N_t}
$.
Inserting the DEIM approximation $\bm{f}(V\bm{\alpha}_{(2)}) \approx W \bm{f}_D(V\bm{\alpha}_{(2)})$ with the reduced nonlinear vector $\bm{f}_D(V\bm{\alpha}_{(2)})=\mathsf{P}^T\bm{f}(V\bm{\alpha}_{(2)})$, into the full discrete P-ROM \eqref{romfullkdv2}, the full discrete PD-ROM for the KdV equation read as
\begin{equation}\label{romfullkdv3}
\delta_t\bm{\alpha}_{(2)}^k+\gamma^2 \widehat{D}^3A_t\bm{\alpha}_{(2)}^k = -\widehat{D}V^TW\int_0^1 \bm{f}_D((1-\xi)V\bm{\alpha}_{(2)}^k + \xi V\bm{\alpha}_{(2)}^{k+1})d\xi .
\end{equation}
The discrete reduced global energy \eqref{vecenergy2} is given as
$$
\widehat{\varepsilon}_h^k = \Delta x \sum_{j=1}^N \left( \frac{\eta}{6}(V\bm{\alpha}_{(2)}^k)_j^3 - \frac{\gamma^2}{2}(DV\bm{\alpha}_{(2)}^k)_j^2\right).
$$

For the numerical tests, we consider the KdV equation on a $P$-periodic interval $[0,P]$ with the spatial mesh size $\Delta x=P/1000$ ($N=1000$)
with parameters  $\gamma=1$ and $\eta=6$.
We compute the initial condition from the exact soliton solution \cite{Eidnes20}
$$
u(x,t)= \frac{1}{2} c \, \mathrm{sech}^2\left(\xi (x,t) - \frac{P}{2}\right), \qquad \xi (x,t) = (-x+c t) \text{ }\mathrm{mod} \text{ } P,
$$
moving with a constant speed $c$ in the positive $x$-direction while keeping its initial shape.
In our numerical experiment, we set $c=4$ and $P=20$. We run the simulation until the final time $T=10$ with the time  step-size $\Delta t=0.01$ ($N_t=1000$).
With these setting, we obtain a snapshot matrix $\mathcal{U}$ of size $1000\times 2000$.

\begin{figure}[htb!]
\centering
\includegraphics[width=0.5\columnwidth]{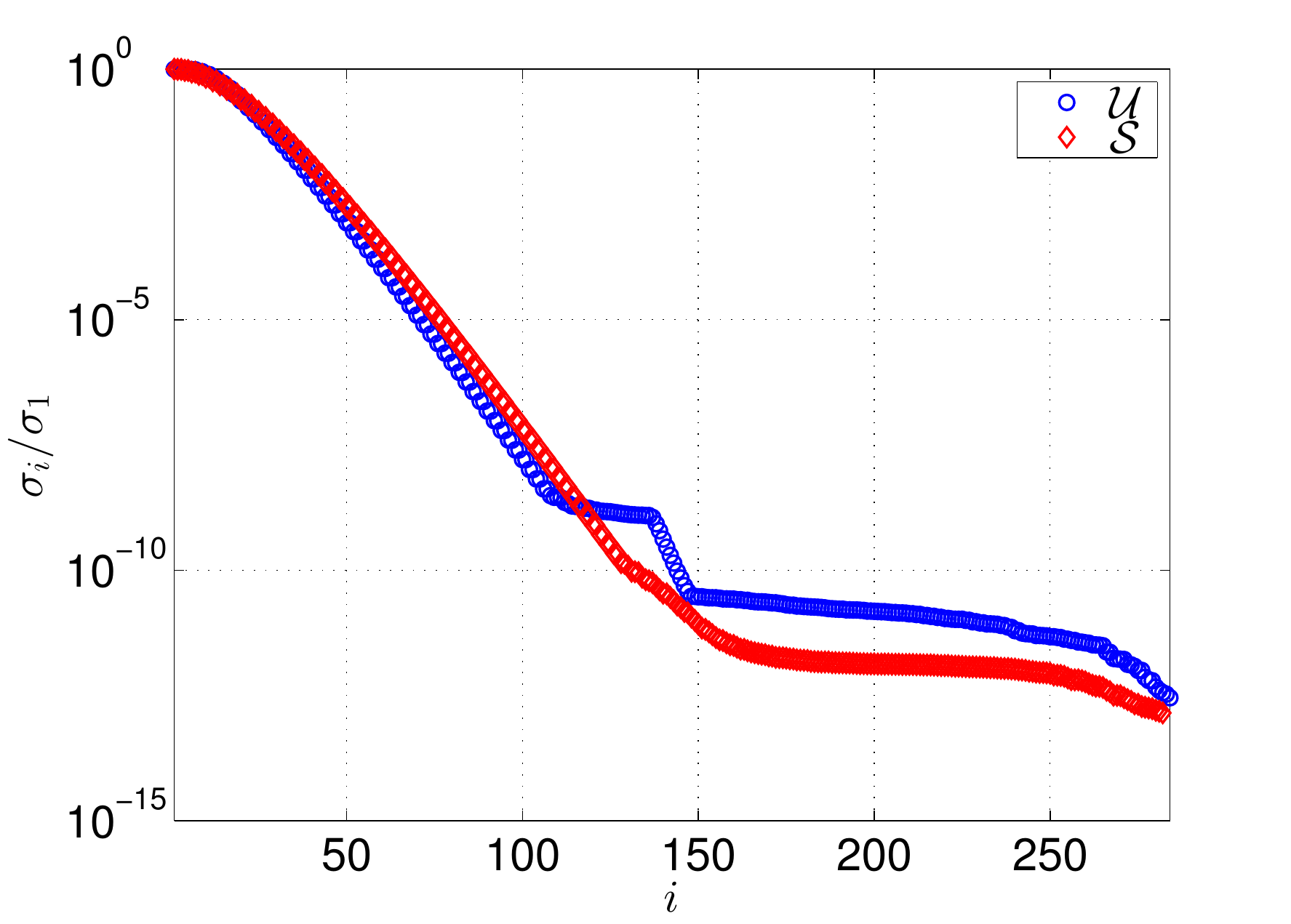}
\caption{KdV: Decay of singular values\label{kdv:sing}}
\end{figure}

The decay of the singular values of the solution snapshot $\mathcal{U}$ and the nonlinear snapshot $\mathcal{S}$ is shown in \figurename~\ref{kdv:sing}.
The slow decay of the singular values is characteristic for the problems with complex wave phenomena.
The Kolmogorov $n$-width of the solution manifolds is a classical concept of nonlinear approximation theory as it describes the error arising from a projection onto the best-possible space of a given dimension $n$.
Hamiltonian systems with soliton solutions have the same characteristics as the wave phenomena and transport dominated problems, i.e., slow decay of the Kolmogorov $n$-widths, e.g., $n^{-1/2}$ \cite{Ohlberger16,Urban19}. We remark that the singular values depend on the snapshots and only provide a bound on the projection error corresponding to the POD space, therefore, they do not, in general, correspond to the Kolmogorov $n$-widths
\cite{Peherstorfer22}.

\begin{figure}[H]
\centering
\includegraphics[width=0.48\columnwidth]{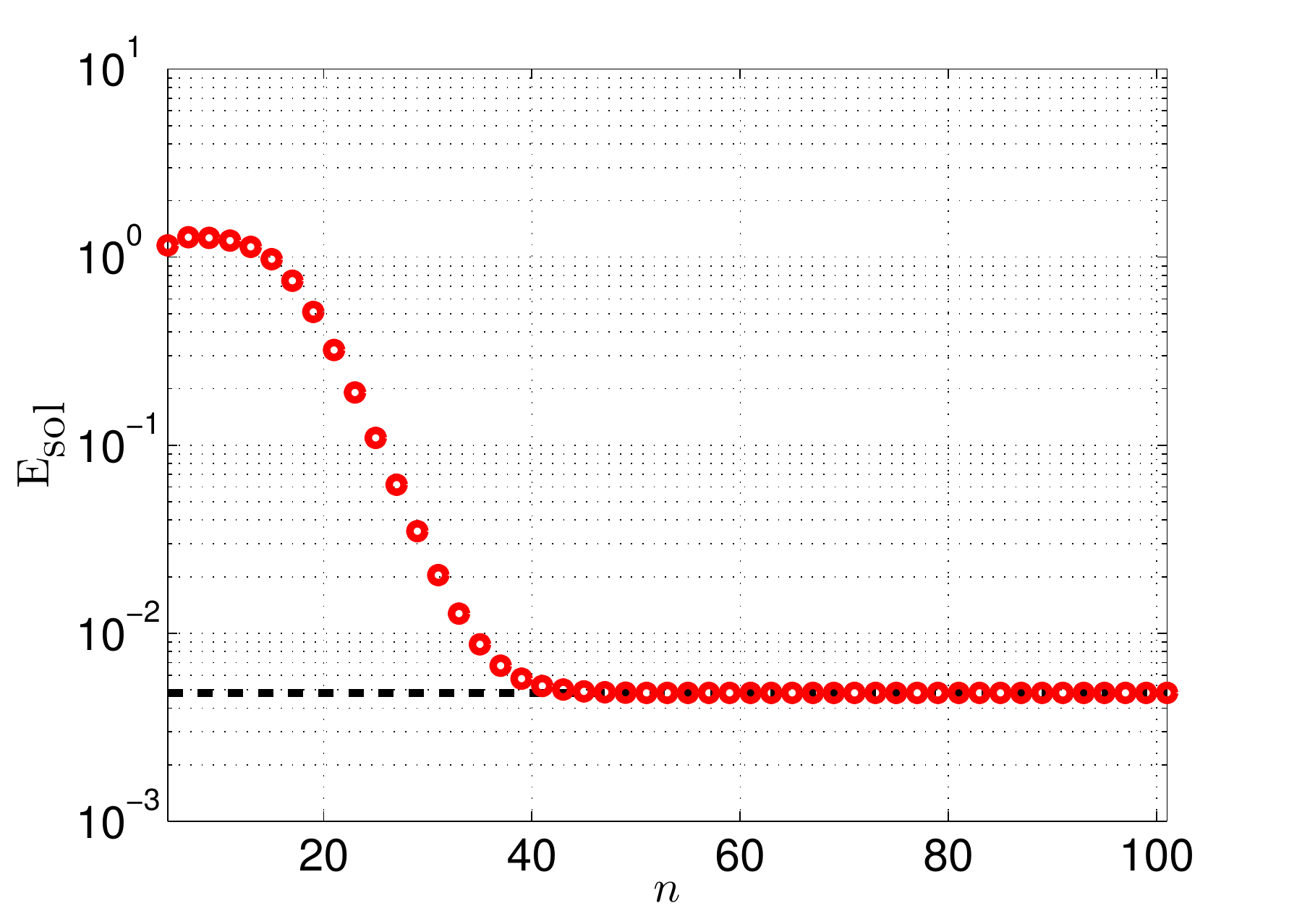}
\includegraphics[width=0.48\columnwidth]{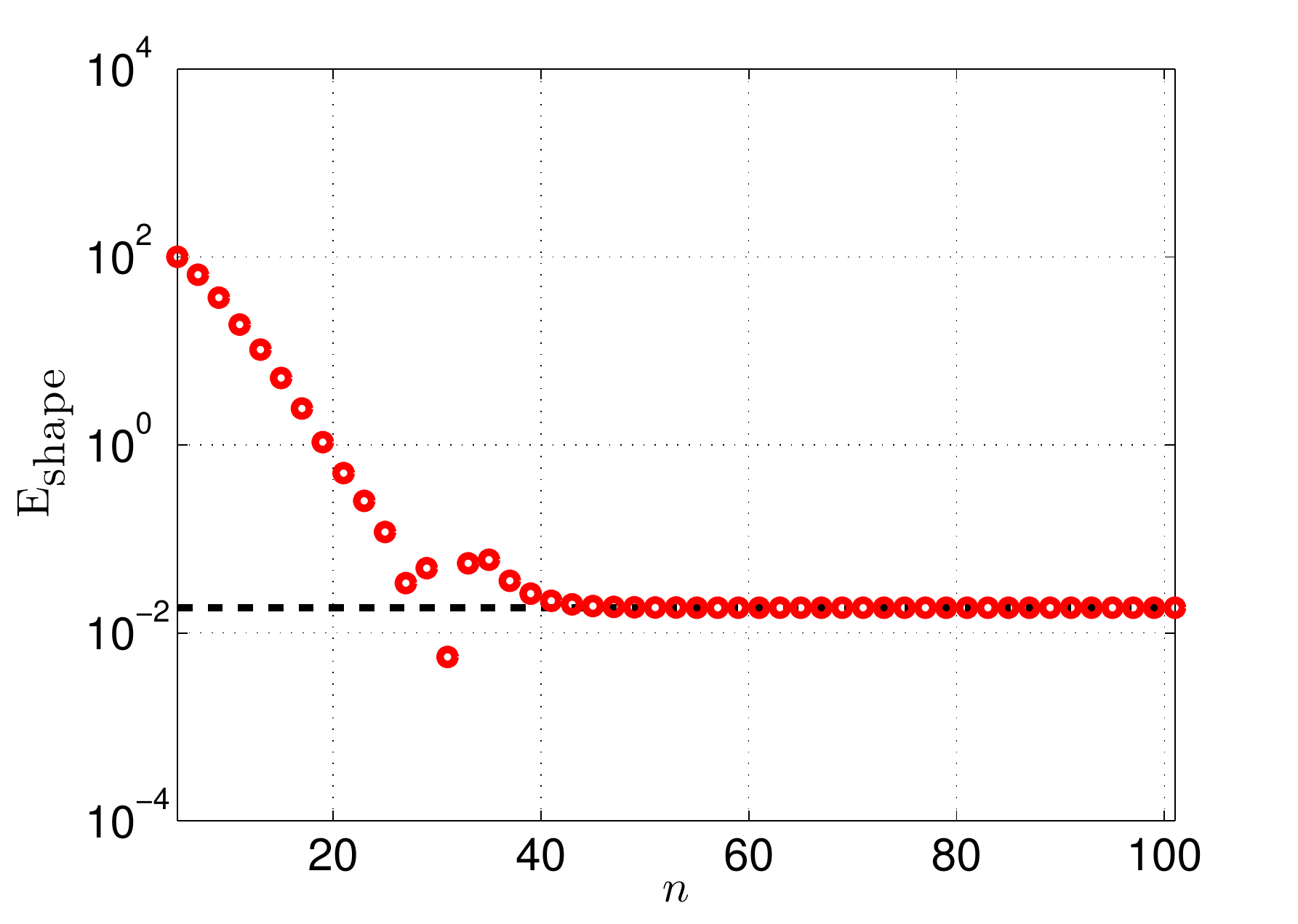}
\caption{KdV: Solution errors $\text{E}_{\text{sol}}$ and shape errors $\text{E}_{\text{shape}}$: dashed line corresponds to the FOM\label{kdv:bench1}}
\end{figure}

\begin{figure}[H]
\centering
\includegraphics[width=0.48\columnwidth]{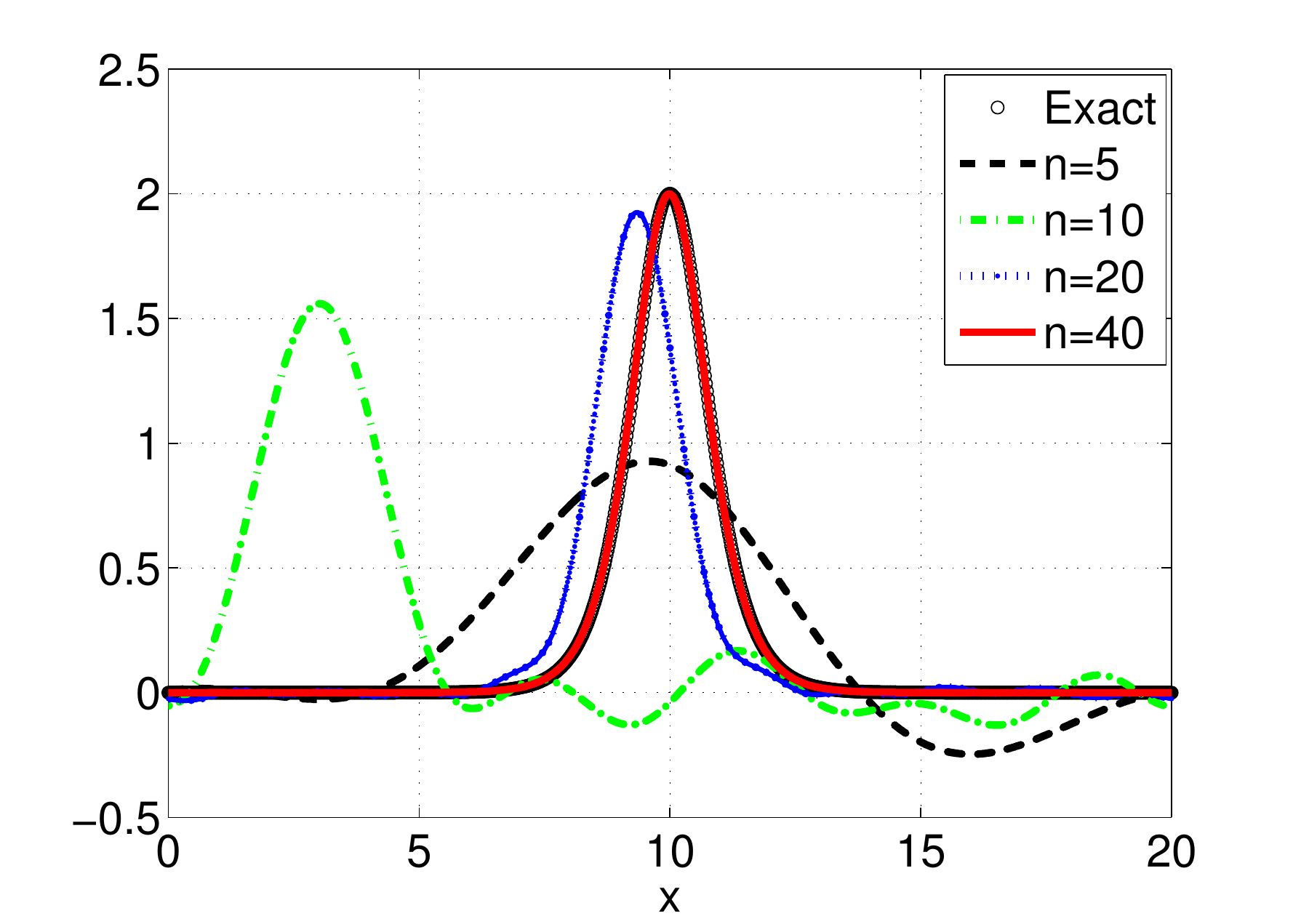}
\includegraphics[width=0.48\columnwidth]{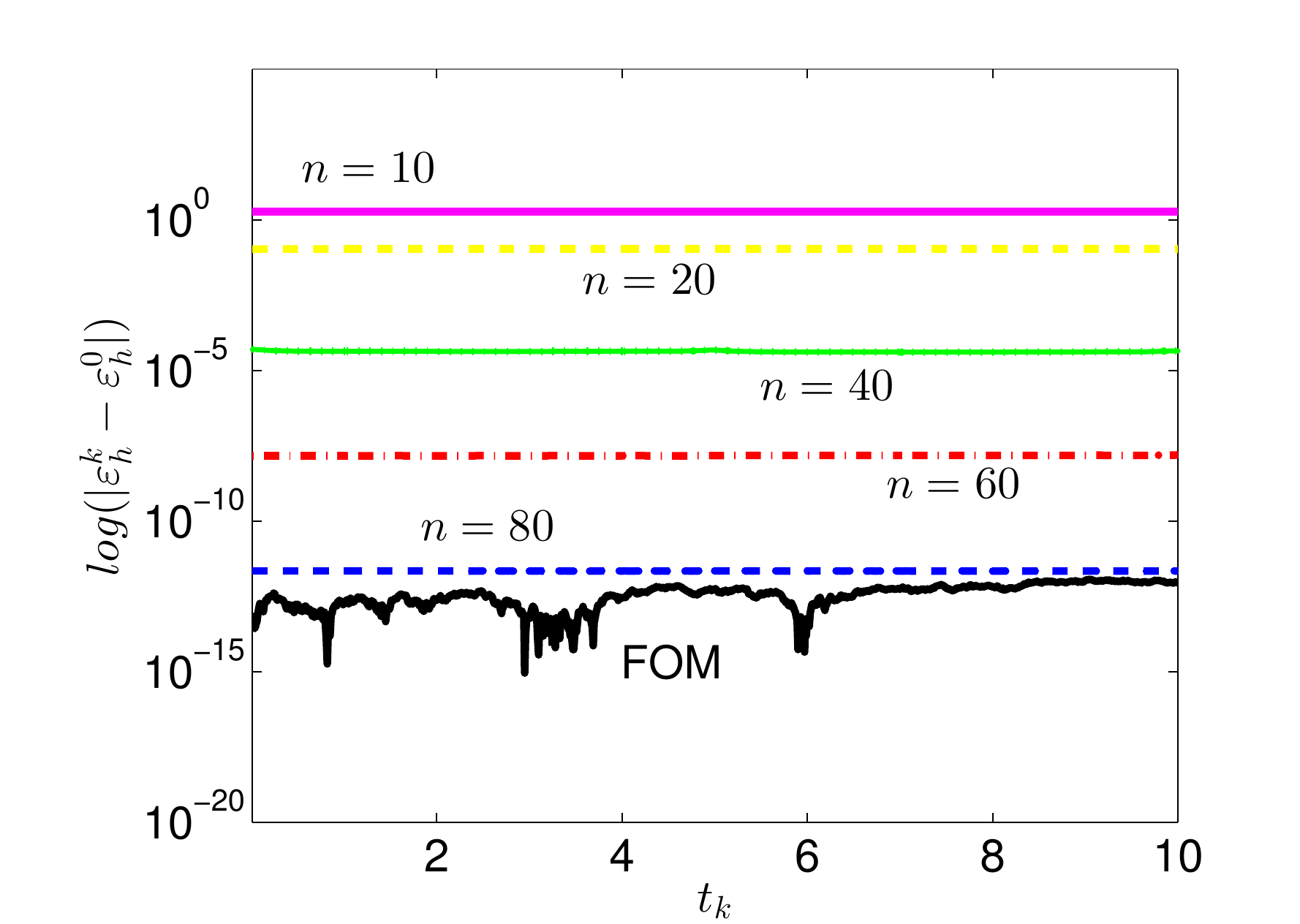}
\caption{KdV: Solitons at the final time and energy preservation for different number of POD modes\label{kdv:bench2}}
\end{figure}

We investigate the accuracy and computational efficiency of the P-ROM \eqref{romfullkdv2} by increasing the number of POD modes up to $n=100$, because the singular values are much smaller than the prescribed tolerances afterward in Figure \ref{kdv:sing}.
In \figurename~\ref{kdv:bench1}, the solution errors $\text{E}_{\text{sol}}$ and the shape errors $\text{E}_{\text{shape}}$ are shown with increasing number of POD modes.
Both errors reach a plateau for $n\ge 40$, therefore we investigate the behavior of the soliton solutions of the ROMs and FOMs for smaller number of POD modes in \figurename~\ref{kdv:bench2}, left, plotted at the final time $t=10$.
We see that for $n=5,10,20$, the exact shape of the solitons cannot be caught by the reduced approximation, while $n=40$ POD modes reproduces almost the exact shape of the FOM soliton.
The discrete reduced global energy in \figurename~\ref{kdv:bench2}, right, is preserved in a small band for either number of POD modes, and they converge to the discrete global energy by increasing number of POD modes.

\begin{figure}[H]
\centerline{\includegraphics[width=1.3\columnwidth]{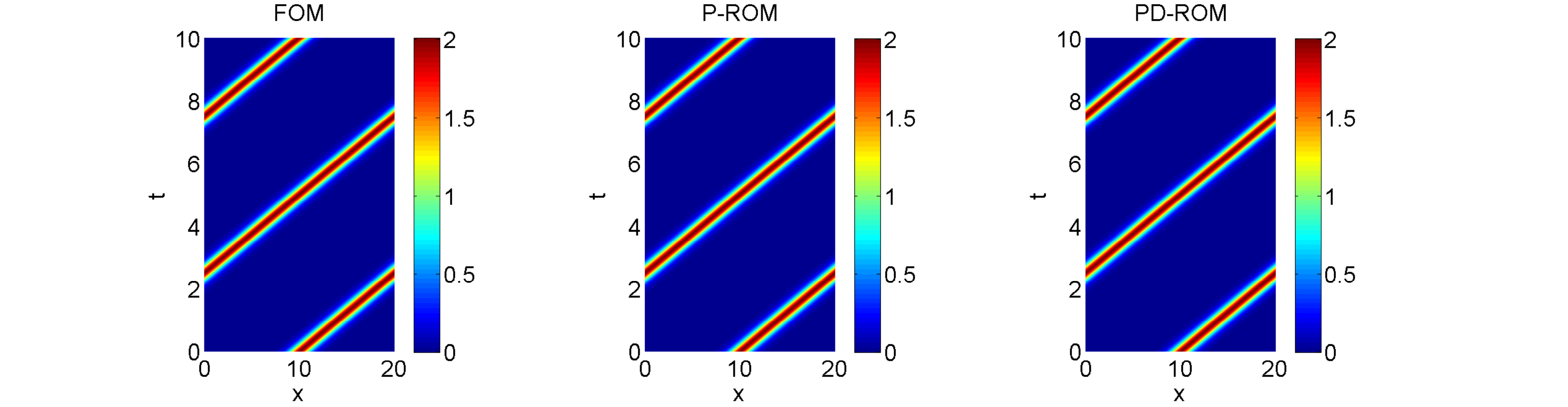}}
\caption{KdV: Time evolution of the soliton waves\label{kdv:surf}}
\end{figure}

\begin{figure}[H]
\centering
\includegraphics[width=0.48\columnwidth]{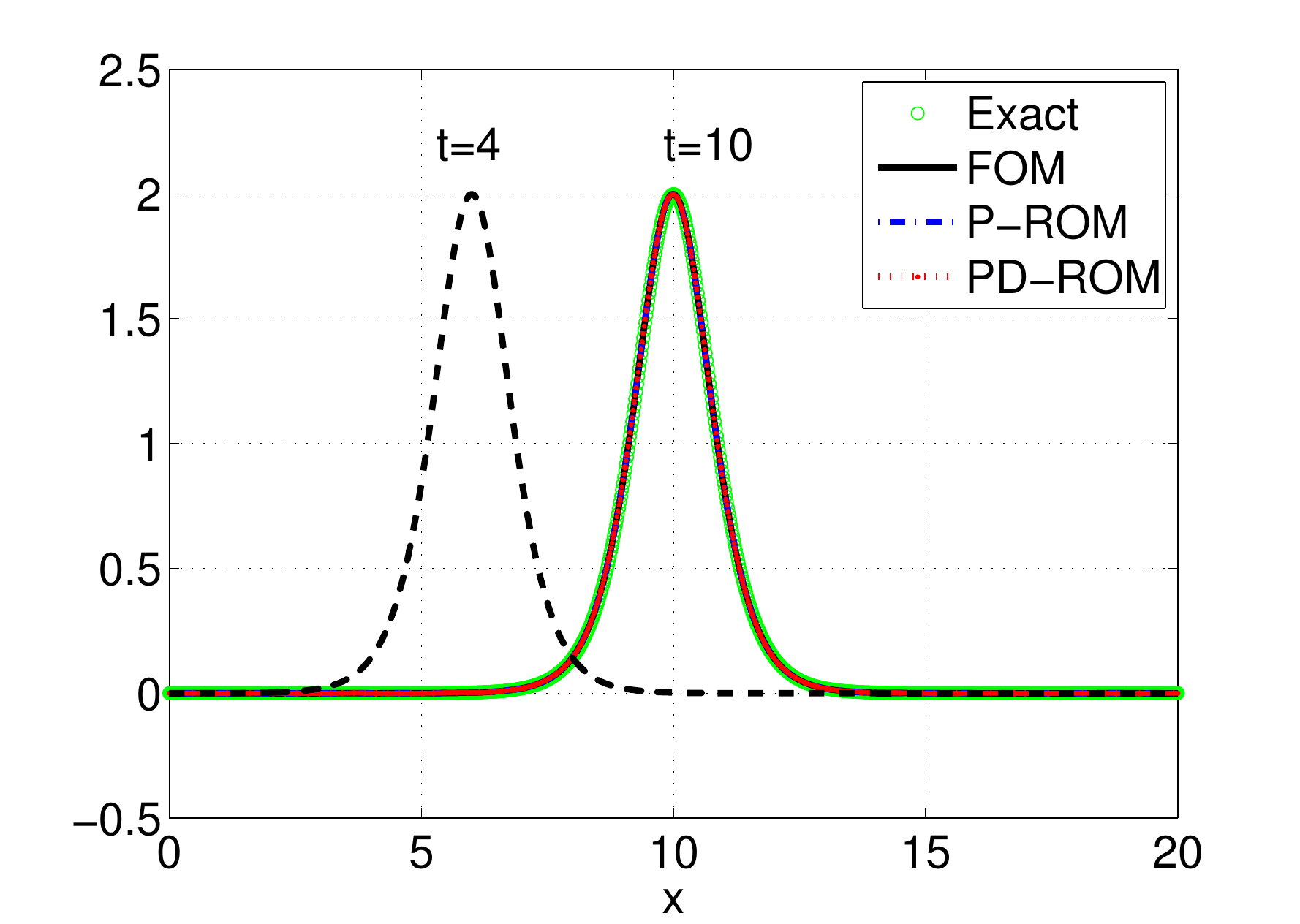}
\includegraphics[width=0.48\columnwidth]{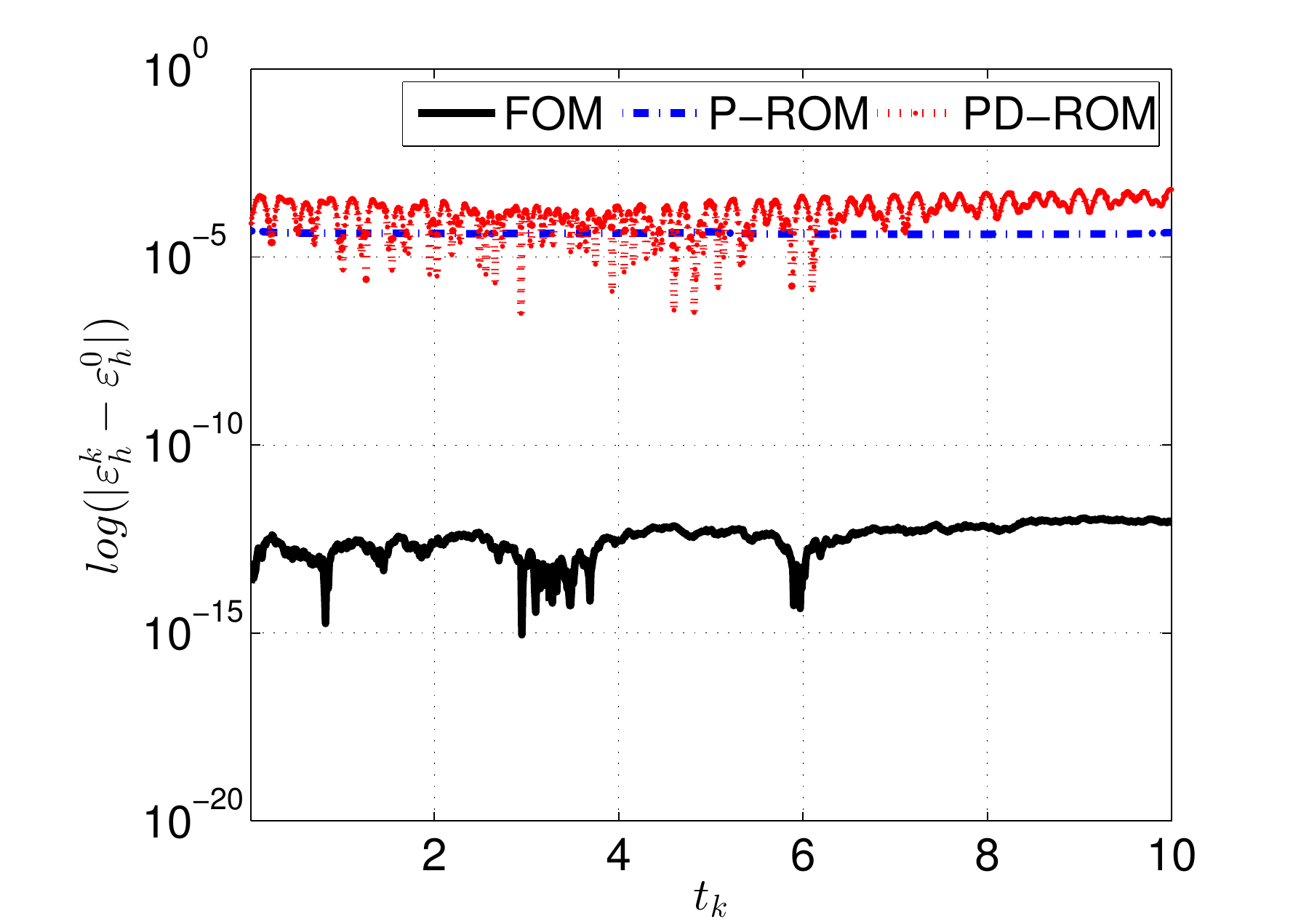}
\caption{KdV: (left) exact solitons at $t=4,10$ and numerical solitons at $t=10$, (right) conservation of discrete global energy\label{kdv:energy}}
\end{figure}

We take the number of POD modes $n=40$ and the number of DEIM modes $\tilde{n}=45$ using the prescribed tolerances.
The solitons of the FOMs and ROMs  behave similar in \figurename~\ref{kdv:surf}.
The exact and numerical solitons at time $t=4$ and at the final time $t=10$ in \figurename~\ref{kdv:energy}, left, show that the shape of the solitons are well preserved.
In \figurename~\ref{kdv:energy}, right, we show the discrete global energy conservation for FOM and ROMs.
The discrete reduced global energy is conserved with increasing number of POD modes, and oscillates in a small band.
This validates the preservation of the discrete reduced global energy in the full-rank reduce space as stated in Theorem~\ref{thm1}.

%%%%%%%%%%%%%%%%%%%%%%%%%%%%%%%%%%%%%%%%%%%%%%%%%%%%%%%%%%%%%%%%%%%%%%%%%%%%%%%%%%%%%%%%%%%%
\subsection{1D NLS equation}
\label{ex:nls}

As the second test problem, we consider one-dimensional NLS equation \eqref{nls1} in multi-symplectic form \eqref{nlsms1}.
The full discrete FOM \eqref{msfull2} for the multi-symplectic NLS equation \eqref{nlsms1} reads as
\begin{equation}\label{fullnlsms}
\begin{aligned}
\delta_t\bm{q}^k - DA_t\bm{v}^k &= \int_0^1 \bm{f}((1-\xi)\bm{p}^k + \xi\bm{p}^{k+1},(1-\xi)\bm{q}^k + \xi\bm{q}^{k+1})d\xi,\\
-\delta_t\bm{p}^k- DA_t\bm{w}^k &= \int_0^1 \bm{g}((1-\xi)\bm{p}^k + \xi\bm{p}^{k+1},(1-\xi)\bm{q}^k + \xi\bm{q}^{k+1})d\xi ,\\
DA_t\bm{p}^k &= A_t\bm{v}^k,\\
DA_t\bm{q}^k &= A_t\bm{w}^k,
\end{aligned}
\end{equation}
where $\bm{p}=\bm{z}_{(1)}$, $\bm{q}=\bm{z}_{(2)}$, $\bm{v}=\bm{z}_{(3)}$ and $\bm{w}=\bm{z}_{(4)}$ are the discrete solutions, $\bm{f}(\bm{p},\bm{q})=\beta (\bm{p}^2+\bm{q}^2)\bm{p}$ and $\bm{g}(\bm{p},\bm{q})=\beta (\bm{p}^2+\bm{q}^2)\bm{q}$ are the nonlinear terms. By elimination of the auxiliary variables $\bm{v}$ and $\bm{w}$, the system \eqref{fullnlsms} can be solved for the unknowns $\bm{p}$ and $\bm{q}$
\begin{equation}\label{fullnlsms2}
\begin{aligned}
\delta_t\bm{q}^k - D^2A_t\bm{p}^k &= \int_0^1 \bm{f}((1-\xi)\bm{p}^k + \xi\bm{p}^{k+1},(1-\xi)\bm{q}^k + \xi\bm{q}^{k+1})d\xi,\\
-\delta_t\bm{p}^k- D^2A_t\bm{q}^k &= \int_0^1 \bm{g}((1-\xi)\bm{p}^k + \xi\bm{p}^{k+1},(1-\xi)\bm{q}^k + \xi\bm{q}^{k+1})d\xi .
\end{aligned}
\end{equation}
The discrete global energy \eqref{vecenergy} is
$$
\varepsilon_h^k = \frac{\Delta x}{2} \sum_{j=1}^N \left( \frac{\beta}{2}\left((\bm{p}_j^k)^2+(\bm{q}_j^k)^2\right)^2 - (D\bm{p}^k)_j^2  - (D\bm{q}^k)_j^2 \right).
$$

The ROMs are computed similar to the KdV equation by eliminating the auxiliary variables  and solving a coupled system for only the states $\bm{p}$ and $\bm{q}$. For the states $p$ and $v$, we take the POD basis matrix  $V_{(1)}\in\mathbb{R}^{N\times n}$, and we take the POD basis matrix  $V_{(2)}\in\mathbb{R}^{N\times n}$ for the states $q$ and $w$
$$
\bm{p}(t) \approx V_{(1)}\bm{\alpha}_{(1)}(t) , \quad \bm{q}(t)\approx V_{(2)}\bm{\alpha}_{(2)}(t), \quad \bm{v}(t)\approx V_{(1)}\bm{\alpha}_{(3)}(t), \quad \bm{w}(t)\approx V_{(2)}\bm{\alpha}_{(4)}(t),
$$
with the reduced coefficient vectors $\bm{\alpha}_{(i)}(t): [0,T]\mapsto \mathbb{R}^{n}$, $i=1,\ldots ,4$.
The snapshot matrices $\mathcal{P}$ and $\mathcal{Q}$ are
\begin{align*}
\mathcal{P} &= [\bm{p}^1, \; \ldots \; ,\bm{p}^{N_t} \; D\bm{p}^1,\; \ldots \;, D\bm{p}^{N_t} ] \in\mathbb{R}^{N\times 2N_t}, \\
\mathcal{Q} &= [\bm{q}^1, \; \ldots \; ,\bm{q}^{N_t} \; D\bm{q}^1, \; \ldots \;,  D\bm{q}^{N_t} ] \in\mathbb{R}^{N\times 2N_t}.
\end{align*}
After elimination of the coefficient vectors $\bm{\alpha}_{(3)}$ and $\bm{\alpha}_{(4)}$ of the auxiliary states $\bm{v}$ and $\bm{w}$, the full discrete P-ROM depends only on  $\bm{\alpha}_{(1)}$ and $\bm{\alpha}_{(2)}$
\begin{equation}\label{romfullnls2}
\begin{aligned}
V_{(1)}^T &V_{(2)}\delta_t\bm{\alpha}_{(2)}^k - \widehat{D}_{(1)}^2A_t\bm{\alpha}_{(1)}^k = \\
&V_{(1)}^T \int_0^1 \bm{f}((1-\xi)V_{(1)}\bm{\alpha}_{(1)}^k + \xi V_{(1)}\bm{\alpha}_{(1)}^{k+1},(1-\xi)V_{(2)}\bm{\alpha}_{(2)}^k + \xi V_{(2)}\bm{\alpha}_{(2)}^{k+1})d\xi,\\
-V_{(2)}^T &V_{(1)}\delta_t\bm{\alpha}_{(1)}^k - \widehat{D}_{(2)}^2A_t\bm{\alpha}_{(2)}^k =\\
& V_{(2)}^T \int_0^1 \bm{g}((1-\xi)V_{(1)}\bm{\alpha}_{(1)}^k + \xi V_{(1)}\bm{\alpha}_{(1)}^{k+1},(1-\xi)V_{(2)}\bm{\alpha}_{(2)}^k + \xi V_{(2)}\bm{\alpha}_{(2)}^{k+1})d\xi,
\end{aligned}
\end{equation}
where $\widehat{D}_{(i)}=V_{(i)}^TDV_{(i)}$, $i=1,2$, are the reduced skew-symmetric matrices.

The PD-ROM for the NLS equation is obtained by approximating the nonlinear vectors $\bm{f}(\bm{p},\bm{q})=\beta (\bm{p}^2+\bm{q}^2)\bm{p}$ and $\bm{g}(\bm{p},\bm{q})=\beta (\bm{p}^2+\bm{q}^2)\bm{q}$ from the column space of the nonlinear snapshot matrices, respectively,
$$
\mathcal{F} = [\bm{f}^1 \; \ldots \; \bm{f}^{N_t}  ] \in\mathbb{R}^{N\times N_t},
\qquad
\mathcal{G} = [\bm{g}^1 \; \ldots \; \bm{g}^{N_t}  ] \in\mathbb{R}^{N\times N_t}.
$$
Inserting the DEIM approximations $\bm{f} \approx W_f \bm{f}_D$ and $\bm{g} \approx W_g \bm{g}_D$ with the reduced nonlinear vectors $\bm{f}_D=\mathsf{P}_f^T\bm{f}$ and $\bm{g}_D=\mathsf{P}_g^T\bm{g}$, into the full discrete P-ROM \eqref{romfullnls2}, the full discrete PD-ROM reads as
\begin{equation}\label{romfullnls3}
\begin{aligned}
V_{(1)}^T&V_{(2)}\delta_t\bm{\alpha}_{(2)}^k - \widehat{D}_{(1)}^2A_t\bm{\alpha}_{(1)}^k = \\
&V_{(1)}^T W_f \int_0^1 \bm{f}_D((1-\xi)V_{(1)}\bm{\alpha}_{(1)}^k + \xi V_{(1)}\bm{\alpha}_{(1)}^{k+1},(1-\xi)V_{(2)}\bm{\alpha}_{(2)}^k + \xi V_{(2)}\bm{\alpha}_{(2)}^{k+1})d\xi,\\
-V_{(2)}^T&V_{(1)}\delta_t\bm{\alpha}_{(1)}^k - \widehat{D}_{(2)}^2A_t\bm{\alpha}_{(2)}^k =\\
& V_{(2)}^T W_g \int_0^1 \bm{g}_D((1-\xi)V_{(1)}\bm{\alpha}_{(1)}^k + \xi V_{(1)}\bm{\alpha}_{(1)}^{k+1},(1-\xi)V_{(2)}\bm{\alpha}_{(2)}^k + \xi V_{(2)}\bm{\alpha}_{(2)}^{k+1})d\xi.
\end{aligned}
\end{equation}

The discrete reduced global energy \eqref{vecenergy2} yields
$$
\widehat{\varepsilon}_h^k = \frac{\Delta x}{2} \sum_{j=1}^N \left( \frac{\beta}{2}\left((V_{(1)}\bm{\alpha}_{(1)}^k)_j^2+(V_{(2)}\bm{\alpha}_{(2)}^k)_j^2\right)^2 - (DV_{(1)}\bm{\alpha}_{(1)}^k)_j^2  - (DV_{(2)}\bm{\alpha}_{(2)}^k)_j^2 \right).
$$

We consider the NLS equation on the space-time domain $[-20,60]\times [0,5]$, and with $\beta=2$.
The initial condition is given by $\psi(x,0) = \sech(x)\exp(ix)$, computed with the analytic solution $\psi(x,t) = \sech(x-2t)\exp(ix)$ \cite{Cai19,Gong14}.
We choose $h=80/1000$ ($N=1000$) and $\Delta t = 0.01$ ($N_t=500$) as spatial and temporal mesh sizes, resulting in snapshot matrices $\mathcal{P}$ and $\mathcal{Q}$ of size $1000\times 1000$.

\begin{figure}[H]
\centering
\includegraphics[width=0.5\columnwidth]{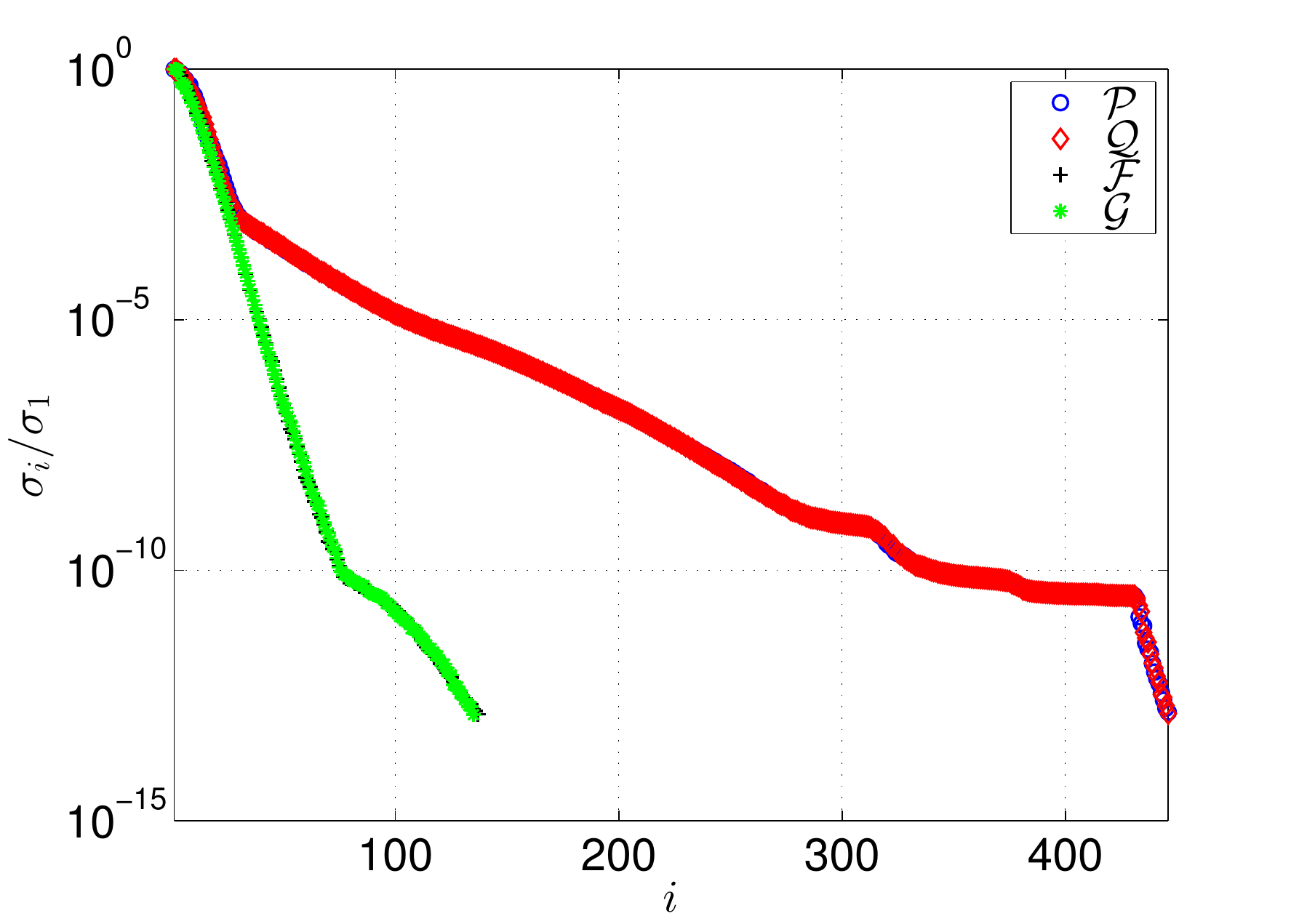}
\caption{1D NLS: Decay of singular values\label{nls1:sing}}
\end{figure}

The decay of the singular values of the solution snapshots $\mathcal{P}$, $\mathcal{Q}$ and the nonlinear snapshots $\mathcal{F}$, $\mathcal{G}$ is given in \figurename~\ref{nls1:sing}.
For the prescribed tolerances, the number of POD and DEIM modes are set as $n=25$ and $\tilde{n}=45$, respectively.

\begin{figure}[H]
\centerline{\includegraphics[width=1.3\columnwidth]{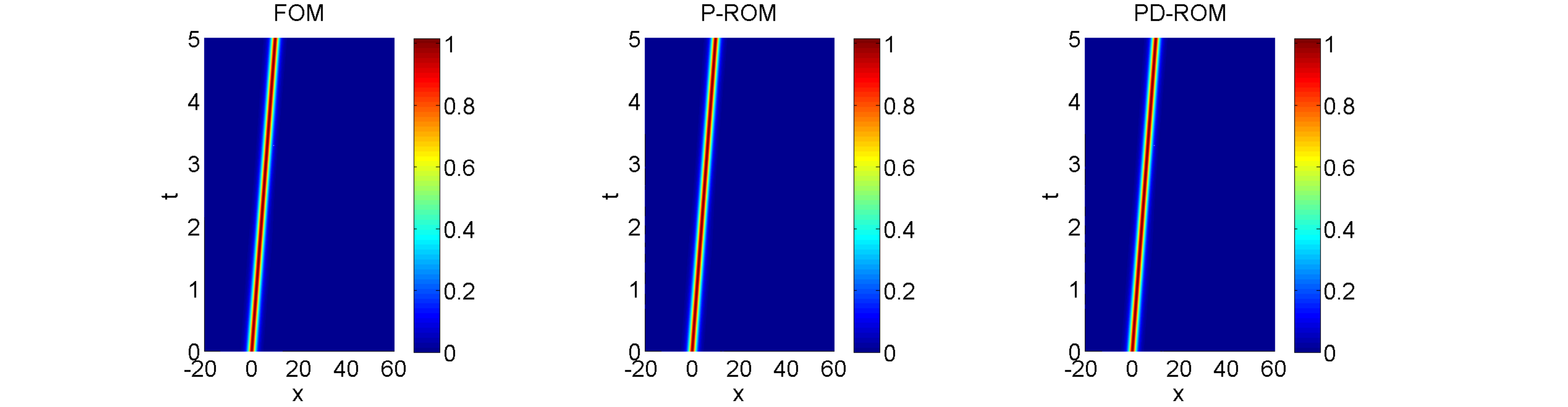}}
	\caption{1D NLS: Time evolution of the soliton waves\label{nls1:surf}}
\end{figure}

\begin{figure}[H]
\centering
\includegraphics[width=0.48\columnwidth]{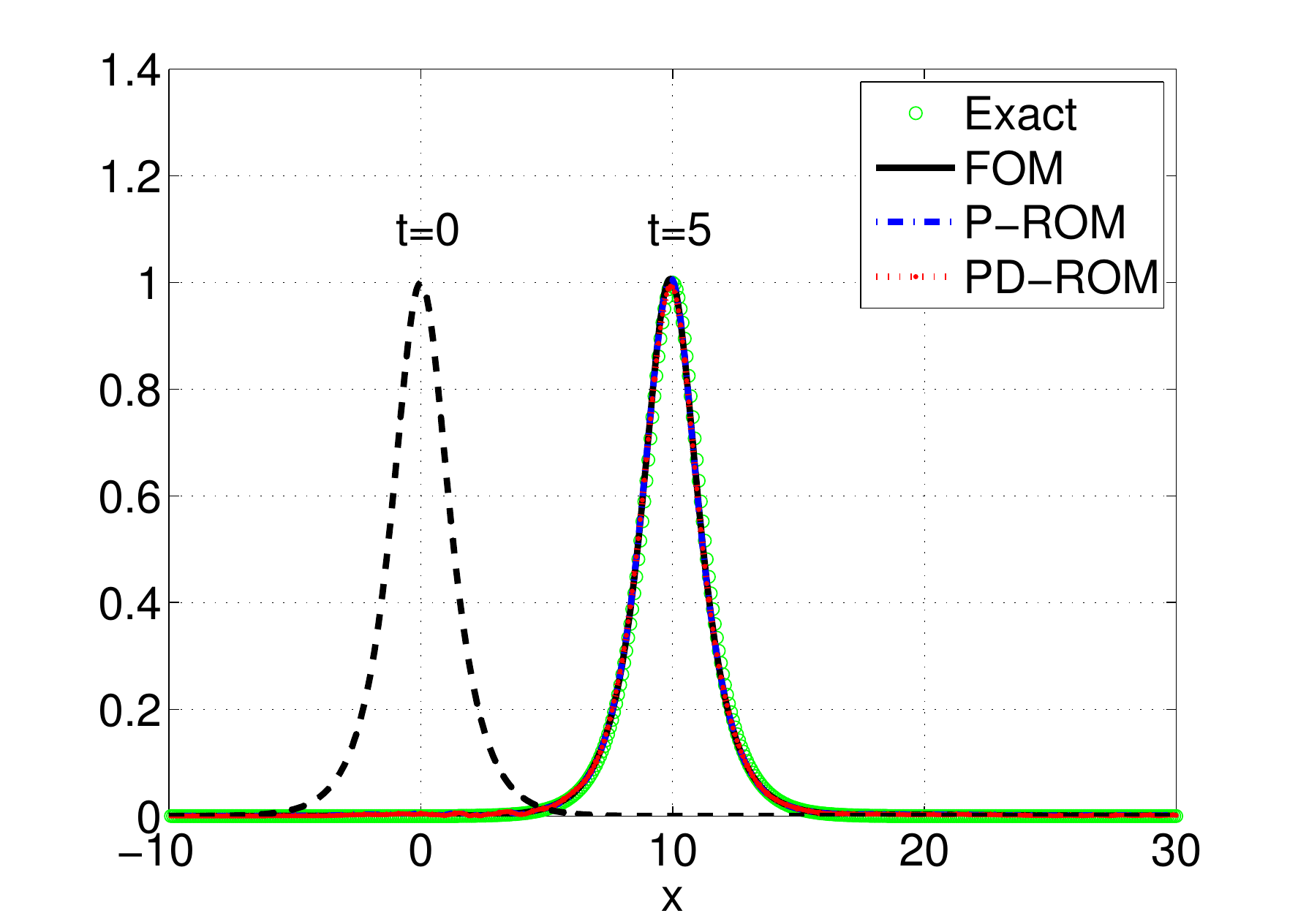}
\includegraphics[width=0.48\columnwidth]{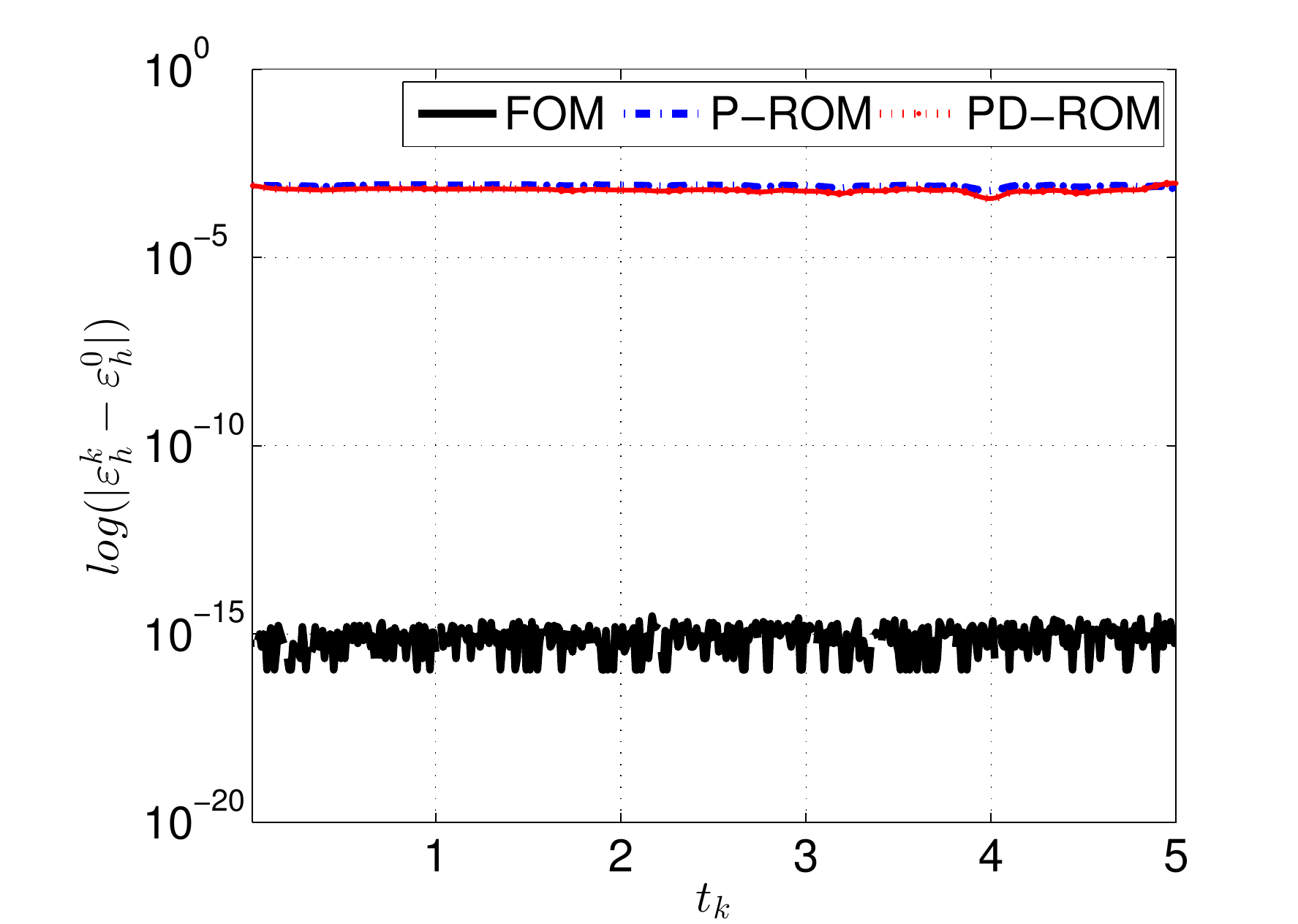}
\caption{1D NLS: (left) exact solitons at $t=0,5$ and numerical solitons at $t=5$, (right) conservation of discrete global energy\label{nls1:energy}}
\end{figure}

Both ROMs behave similar in \figurename~\ref{nls1:surf}.
In \figurename~\ref{nls1:energy}, left, we give the exact and numerical wave plots at the initial time and at the final time.
The shape and the speed of the solitons are well-preserved, see also \tablename~\ref{tbl:error}.
The reduced global energies in \figurename~\ref{nls1:energy}, right, are well preserved with small oscillations.

%%%%%%%%%%%%%%%%%%%%%%%%%%%%%%%%%%%%%%%%%%%%%%%%%%%%%%%%5
\subsection{Zakharov-Kuznetsov equation}

The ZK equation is the two-dimensional generalization of the KdV equation
\begin{equation}\label{zk1}
u_t + uu_{x} + u_{xxx} + u_{xyy}  = 0.
\end{equation}
The ZK equation \eqref{zk1} in multi-symplectic form is given as \cite{Bridges01}
\begin{equation}\label{zkms1}
\begin{aligned}
\phi_x &= u,\\
\frac{1}{2}\phi_t + v_x + w_y &= p - \frac{1}{2}u^2,\\
w_x-v_y&=0,\\
-\frac{1}{2}u_t-p_x&=0,\\
-u_x+q_y &= -v, \\
-q_x-u_y &= -w,
\end{aligned}
\end{equation}
with the Hamiltonian
$S(z)=up-(v^2+w^2)/2-u^3/6$, and the skew-symmetric matrices
\begin{align*}
K&=\begin{bmatrix}
    0 & 0 & 0 & 0 & 0 & 0 \\
    0 & 0 & 0 & 1/2 & 0 & 0 \\
    0 & 0 & 0 & 0 & 0 & 0 \\
    0 & -1/2 & 0 & 0 & 0 & 0 \\
		0 & 0 & 0 & 0 & 0 & 0 \\
		0 & 0 & 0 & 0 & 0 & 0 \\
\end{bmatrix},\\
L_1&=
\begin{bmatrix}
    0 & 0 & 0 & 1 & 0 & 0 \\
    0 & 0 & 0 & 0 & 1 & 0 \\
    0 & 0 & 0 & 0 & 0 & 1 \\
    -1 & 0 & 0 & 0 & 0 & 0 \\
		0 & -1 & 0 & 0 & 0 & 0 \\
		0 & 0 & -1 & 0 & 0 & 0 \\
\end{bmatrix},\qquad
L_2=
\begin{bmatrix}
    0 & 0 & 0 & 0 & 0 & 0 \\
    0 & 0 & 0 & 0 & 0 & 1 \\
    0 & 0 & 0 & 0 & -1 & 0 \\
    0 & 0 & 0 & 0 & 0 & 0 \\
		0 & 0 & 1 & 0 & 0 & 0 \\
		0 & -1 & 0 & 0 & 0 & 0 \\
\end{bmatrix}.
\end{align*}
Under periodic boundary condition, the global energy
\begin{equation}\label{zkenergy1}
\varepsilon (t) = \int_{\Omega} \left(  \frac{1}{2}(u_x^2  + u_y^2)-  \frac{1}{6}u^3 \right)dxdy,
\end{equation}
is conserved.

The FOM and ROMs of ZK equation are computed in the same way as for the 1D KdV equation.
After elimination of the auxiliary variables, the full discrete FOM \eqref{msfull2} of the ZK equation is obtained for the unknown vector $\bm{u}$ only
\begin{equation}\label{fullzkms2}
\delta_t\bm{u}^k + D_x(D_x^2+D_y^2)A_t\bm{u}^k = -D_x\int_0^1 \bm{f}((1-\xi)\bm{u}^k + \xi\bm{u}^{k+1})d\xi ,
\end{equation}
where the differential matrices $D_x$ and $D_y$ are given as in \eqref{2dmat}, and $\bm{f}(\bm{u})=\bm{u}^2/2$ is the nonlinear vector.
The discrete global energy for the ZK equation yields
$$
\varepsilon_h^k = \Delta x\Delta y \sum_{j=1}^{N^2} \left( \frac{1}{2}(D_x\bm{u}^k)_j^2 + \frac{1}{2}(D_y\bm{u}^k)_j^2 - \frac{1}{6}(\bm{u}_j^k)^3 \right).
$$

By taking the same POD basis matrix $V\in\mathbb{R}^{N\times n}$ for each states, and by the elimination of the coefficient vectors of the auxiliary states, we obtain the full discrete P-ROM for the unknown vector $\bm{\alpha}_{(2)}$ only
\begin{equation}\label{romfullzk2}
\delta_t\bm{\alpha}_{(2)}^k+\widehat{D}_x(\widehat{D}_x^2+\widehat{D}_y^2)A_t\bm{\alpha}_{(2)}^k = -\widehat{D}_xV^T\int_0^1 \bm{f}((1-\xi)V\bm{\alpha}_{(2)}^k + \xi V\bm{\alpha}_{(2)}^{k+1})d\xi,
\end{equation}
where $\widehat{D}_x=V^TD_xV\in\mathbb{R}^{n\times n}$ and $\widehat{D}_y=V^TD_yV\in\mathbb{R}^{n\times n}$ are the reduced skew-symmetric matrices.
The POD basis matrix $V$ is computed by the application of the POD method to the snapshot matrix
$$
\mathcal{U} = [\bm{u}^1,  \; \ldots \;, \bm{u}^{N_t} \; D_x\bm{u}^1, \; \ldots \; ,D_x\bm{u}^{N_t} \; D_y\bm{u}^1, \; \ldots \;, D_y\bm{u}^{N_t} ] \in\mathbb{R}^{N\times 3N_t}.
$$
The nonlinear snapshot matrix has the form
$
\mathcal{S} = [\bm{f}^1 \; \ldots \; \bm{f}^{N_t}  ] \in\mathbb{R}^{N\times N_t}
$.
Inserting the DEIM approximation $\bm{f} \approx W \bm{f}_D$ into the full discrete P-ROM \eqref{romfullzk2}, the full discrete PD-ROM is obtained as
\begin{equation}\label{romfullzk3}
\delta_t\bm{\alpha}_{(2)}^k+\widehat{D}_x(\widehat{D}_x^2+\widehat{D}_y^2)A_t\bm{\alpha}_{(2)}^k = -\widehat{D}_xV^TW\int_0^1 \bm{f}_D((1-\xi)V\bm{\alpha}_{(2)}^k + \xi V\bm{\alpha}_{(2)}^{k+1})d\xi.
\end{equation}
The discrete reduced global energy \eqref{vecenergy2} for the ZK equation reads as
$$
\varepsilon_h^k = \Delta x\Delta y \sum_{j=1}^{N^2} \left( \frac{1}{2}(D_xV\bm{\alpha}_{(2)}^k)_j^2 + \frac{1}{2}(D_yV\bm{\alpha}_{(2)}^k)_j^2 - \frac{1}{6}(V\bm{\alpha}_{(2)}^k)_j^3 \right).
$$

For simulations, we consider the ZK equation on a $P$-periodic square domain $[0,P]^2$ with the spatial mesh sizes $\Delta x=\Delta y=P/100$ ($N=10000$).
We compute the initial condition from the exact solution \cite{Chen112d}
$$
u(x,y,t)= 3 c \, \mathrm{sech}^2\left(\frac{\sqrt{c}}{2}\left(\xi (x,t) - \frac{P}{2}\right)\right), \qquad \xi (x,t) = (x-c t) \text{ }\mathrm{mod} \text{ } P,
$$
which is a soliton moving with a constant speed $c$ in the positive $x$-direction while keeping its initial shape. We set $c=1$ and $P=20$.
We run the simulation until the final time $T=5$ with the temporal step-size $\Delta t=0.01$ ($N_t=500$).
The snapshot matrix $\mathcal{U}$ is of size $10000\times 1500$.

\begin{figure}[H]
\centering
\includegraphics[width=0.5\columnwidth]{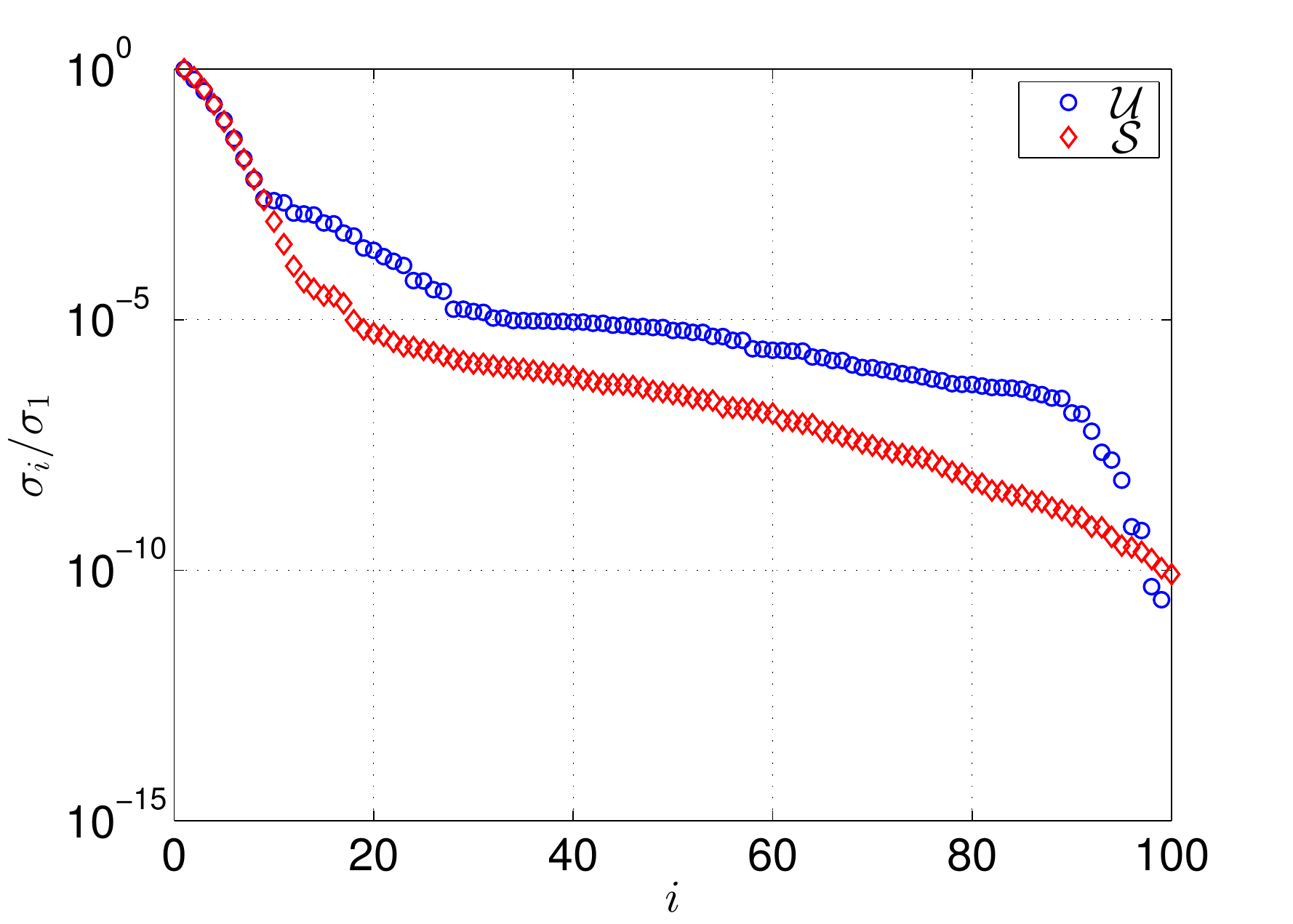}
\caption{ZK: Decay of singular values\label{zk:sing}}
\end{figure}

The decay of the singular values of the solution snapshot $\mathcal{U}$ and the nonlinear snapshot $\mathcal{S}$ is given in \figurename~\ref{zk:sing}.
The number of POD modes and DEIM modes are $n=15$ and $\tilde{n}=25$, respectively.

\begin{figure}[htb!]
\centerline{\includegraphics[width=1.3\columnwidth]{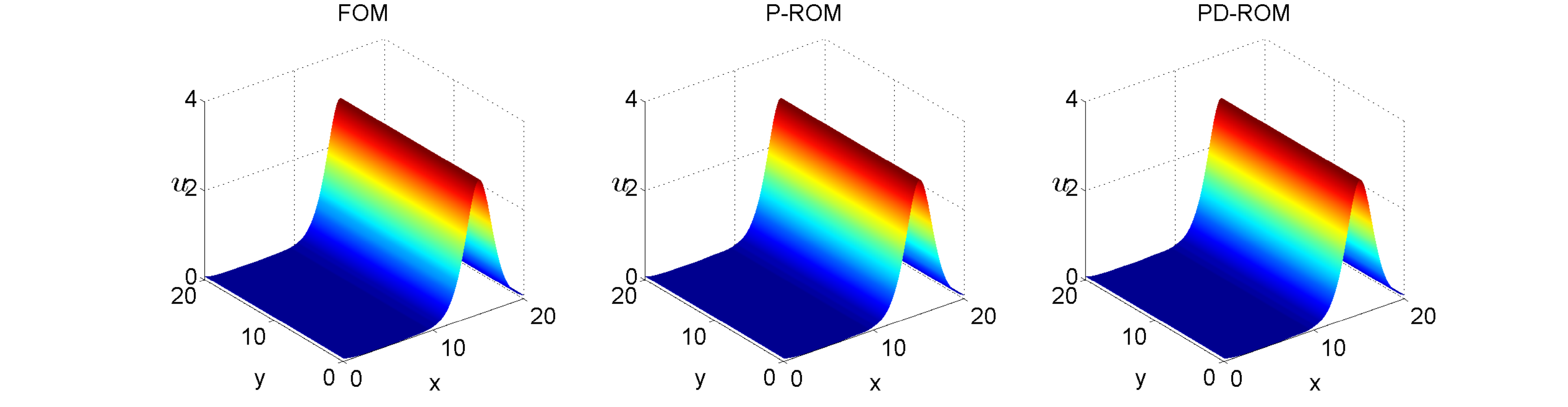}}
\caption{ZK: Soliton waves at $t=5$\label{zk:surf}}
\end{figure}

\begin{figure}[H]
\centering
\includegraphics[width=0.5\columnwidth]{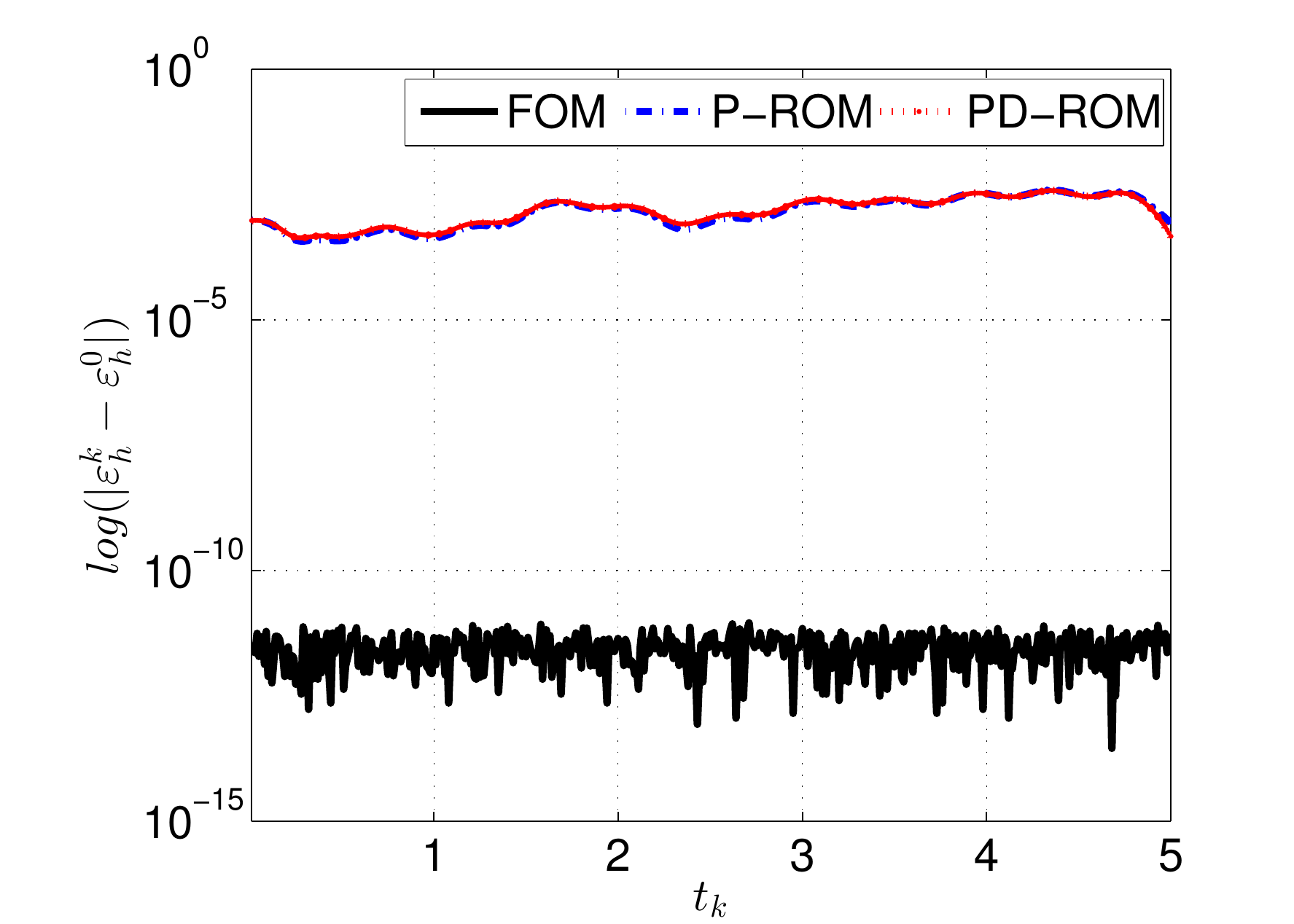}
\caption{ZK: Conservation of discrete global energy\label{zk:energy}}
\end{figure}

The profiles of the solitons of the FOM and ROMs in \figurename~\ref{zk:surf} are very close.
The discrete global energy is conserved by the FOM with a high accuracy, while the discrete reduced global energy by less accurately by P-ROM and PD-ROM in \figurename~\ref{zk:energy}

%%%%%%%%%%%%%%%%%%%%%%%%%%%%%%%%%%%%%%%%%%%%%%%%%%%%%%%%%%%%%%%%%%%%%%%%%%%%%%%%%%%%%%%%%%%%
\subsection{2D NLS equation}

2D NLS equation with two spatial variables is given as
\begin{equation}\label{2dnls1}
i\psi_{t}+\mu(\psi_{xx}+\psi_{yy})+ R'(|\psi|^2,x,y)\psi  =0,
\end{equation}
with the complex-valued solution $\psi(x,y,t)$, where $R'$ denotes the derivative of $R$ with respect to the first argument.
Splitting of the solution into real and imaginary part, $\psi=p+iq$, and introducing the auxiliary variables
$v=\partial_{x}p$, $w=\partial_{x}q$, $a=\partial_{y}p$ and $b=\partial_{y}q$, the 2D NLS equation \eqref{2dnls1} can be written equivalently as \cite{Chen112dnls,Li15nls}
\begin{equation}\label{2dnlsms1}
\begin{aligned}
q_t - \mu v_x - \mu a_y &= R'(p^2+q^2,x,y)p,\\
-p_t - \mu w_x - \mu b_y &= R'(p^2+q^2,x,y)q,\\
\mu p_x &= \mu v,\\
\mu q_x &= \mu w,\\
\mu p_y &= \mu a,\\
\mu q_y &= \mu b,
\end{aligned}
\end{equation}
which is in multi-symplectic form \eqref{ms1} with the Hamiltonian
$S(z)=R(p^{2}+q^{2},x,y)/2+\mu(v^{2}+w^{2}+a^{2}+b^{2})/2$, and the skew-symmetric matrices
\begin{align*}
K &=\left(\begin{array}{cccccc}0&1&0&0&0&0\\-1&0&0&0&0&0\\0&0&0&0&0&0\\0&0&0&0&0&0\\0&0&0&0&0&0\\0&0&0&0&0&0\\\end{array}\right),\\
L_1 &=\left(\begin{array}{cccccc}0&0&-\mu&0&0&0\\0&0&0&-\mu&0&0\\ \mu&0&0&0&0&0\\0&\mu&0&0&0&0\\0&0&0&0&0&0\\0&0&0&0&0&0\\\end{array}\right),
\quad
L_2=\left(\begin{array}{cccccc}0&0&0&0&-\mu&0\\0&0&0&0&0&-\mu\\0&0&0&0&0&0\\0&0&0&0&0&0\\ \mu&0&0&0&0&0\\0&\mu&0&0&0&0\\\end{array}\right).
\end{align*}

Under periodic boundary condition, the system \eqref{2dnlsms1} conserves the global energy
\begin{equation}\label{2dnlsenergy1}
\varepsilon (t) = \frac{1}{2}\int_{\Omega} \left( R(p^2+q^2,x,y) - \mu (p_x^2  + q_x^2+p_y^2  + q_y^2) \right)dxdy.
\end{equation}

The FOM and ROMs of the 2D NLS equation are constructed as the same as the 1D NLS equation.
Eliminating the auxiliary variables in the full discrete formulation \eqref{msfull2}, we  obtain the full discrete FOM with the  unknowns $\bm{p}$ and $\bm{q}$ only
\begin{equation}\label{full2dnlsms2}
\begin{aligned}
\delta_t\bm{q}^k - \mu (D_x^2 +D_y^2)A_t\bm{p}^k  &= \int_0^1 \bm{f}((1-\xi)\bm{p}^k + \xi\bm{p}^{k+1},(1-\xi)\bm{q}^k + \xi\bm{q}^{k+1})d\xi,\\
-\delta_t\bm{p}^k - \mu (D_x^2 +D_y^2)A_t\bm{q}^k &= \int_0^1 \bm{g}((1-\xi)\bm{p}^k + \xi\bm{p}^{k+1},(1-\xi)\bm{q}^k + \xi\bm{q}^{k+1})d\xi,
\end{aligned}
\end{equation}
where the nonlinear vectors are given by $\bm{f}(\bm{p},\bm{q})=\bm{R}'(\bm{p}^2+\bm{q}^2,\bm{x})\bm{p}$ and $\bm{g}(\bm{p},\bm{q})=\bm{R}'(\bm{p}^2+\bm{q}^2,\bm{x})\bm{q}$.
The discrete global energy for the 2D NLS equation yields
$$
\varepsilon_h^k = \frac{\Delta x \Delta y}{2} \sum_{j=1}^N \left[ \bm{R}\left((\bm{p}_j^k)^2+(\bm{q}_j^k)^2,\bm{x}_j\right) - \mu \left( (D_x\bm{p}^k)_j^2  + (D_x\bm{q}^k)_j^2 + (D_y\bm{p}^k)_j^2 + (D_y\bm{q}^k)_j^2 \right)\right].
$$

By taking the same POD basis matrix $V_{(1)}\in\mathbb{R}^{N\times n}$ for the states  $p$, $v$ and $a$, and the same POD basis matrix $V_{(2)}\in\mathbb{R}^{N\times n}$ for the states $q$, $w$ and $b$, and eliminating the coefficient vectors of the auxiliary states, we obtain the full discrete P-ROM with the unknown vectors $\bm{\alpha}_{(1)}$ and $\bm{\alpha}_{(2)}$ only
\begin{equation}\label{romfull2dnls2}
\begin{aligned}
V_{(1)}^T&V_{(2)}\delta_t\bm{\alpha}_{(2)}^k -  \mu (\widehat{D}_{x,(1)}^2+\widehat{D}_{y,(1)}^2)A_t\bm{\alpha}_{(1)}^k =\\
& V_{(1)}^T \int_0^1 \bm{f}((1-\xi)V_{(1)}\bm{\alpha}_{(1)}^k + \xi V_{(1)}\bm{\alpha}_{(1)}^{k+1},(1-\xi)V_{(2)}\bm{\alpha}_{(2)}^k + \xi V_{(2)}\bm{\alpha}_{(2)}^{k+1})d\xi,\\
-V_{(2)}^T&V_{(1)}\delta_t\bm{\alpha}_{(1)}^k -  \mu (\widehat{D}_{x,(2)}^2+\widehat{D}_{y,(2)}^2)A_t\bm{\alpha}_{(2)}^k= \\
& V_{(2)}^T \int_0^1 \bm{g}((1-\xi)V_{(1)}\bm{\alpha}_{(1)}^k + \xi V_{(1)}\bm{\alpha}_{(1)}^{k+1},(1-\xi)V_{(2)}\bm{\alpha}_{(2)}^k + \xi V_{(2)}\bm{\alpha}_{(2)}^{k+1})d\xi,
\end{aligned}
\end{equation}
where $\widehat{D}_{x,(i)}=V_{(i)}^TD_xV_{(i)}$ and $\widehat{D}_{y,(i)}=V_{(i)}^TD_yV_{(i)}$, $i=1,2$, are the reduced skew-symmetric matrices.
The POD basis matrices $V_{(1)}$ and $V_{(2)}$ are computed with the snapshot matrices $\mathcal{P}$ and $\mathcal{Q}$, respectively,
\begin{align*}
\mathcal{P} &= [\bm{p}^1 \; \ldots \; \bm{p}^{N_t} \; D_x\bm{p}^1\; \ldots \; D_x\bm{p}^{N_t}\; D_y\bm{p}^1\; \ldots \; D_y\bm{p}^{N_t} ] \in\mathbb{R}^{N\times 3N_t}, \\
\mathcal{Q} &= [\bm{q}^1 \; \ldots \; \bm{q}^{N_t} \; D_x\bm{q}^1\; \ldots \; D_x\bm{q}^{N_t} \; D_y\bm{q}^1\; \ldots \; D_y\bm{q}^{N_t} ] \in\mathbb{R}^{N\times 3N_t}.
\end{align*}

The PD-ROM for the 2D NLS equation is again obtained by approximating the nonlinear vectors $\bm{f}$ and $\bm{g}$ from the column space of the nonlinear snapshot matrices
$$
\mathcal{F} = [\bm{f}^1 \; \ldots \; \bm{f}^{N_t}  ] \in\mathbb{R}^{N\times N_t} , \qquad \mathcal{G} = [\bm{g}^1 \; \ldots \; \bm{g}^{N_t}  ] \in\mathbb{R}^{N\times N_t}.
$$
Inserting the DEIM approximations $\bm{f} \approx W_f \bm{f}_D$ and $\bm{g} \approx W_g \bm{g}_D$ into the full discrete P-ROM \eqref{romfull2dnls2}, the full discrete PD-ROM for the 2D NLS equation read as
\begin{equation}\label{romfull2dnls3}
\begin{aligned}
V_{(1)}^T&V_{(2)}\delta_t\bm{\alpha}_{(2)}^k -  \mu (\widehat{D}_{x,(1)}^2+\widehat{D}_{y,(1)}^2)A_t\bm{\alpha}_{(1)}^k =\\
& V_{(1)}^T W_f\int_0^1 \bm{f}_D((1-\xi)V_{(1)}\bm{\alpha}_{(1)}^k + \xi V_{(1)}\bm{\alpha}_{(1)}^{k+1},(1-\xi)V_{(2)}\bm{\alpha}_{(2)}^k + \xi V_{(2)}\bm{\alpha}_{(2)}^{k+1})d\xi,\\
-V_{(2)}^T&V_{(1)}\delta_t\bm{\alpha}_{(1)}^k -  \mu (\widehat{D}_{x,(2)}^2+\widehat{D}_{y,(2)}^2)A_t\bm{\alpha}_{(2)}^k =\\
& V_{(2)}^T W_g\int_0^1 \bm{g}_D((1-\xi)V_{(1)}\bm{\alpha}_{(1)}^k + \xi V_{(1)}\bm{\alpha}_{(1)}^{k+1},(1-\xi)V_{(2)}\bm{\alpha}_{(2)}^k + \xi V_{(2)}\bm{\alpha}_{(2)}^{k+1})d\xi,
\end{aligned}
\end{equation}
The discrete reduced global energy \eqref{vecenergy2} is computed as
\begin{align*}
\widehat{\varepsilon}_h^k = \frac{\Delta x \Delta y}{2} &\sum_{j=1}^N [ \bm{R}\left((V_{(1)}\bm{p}^k)_j^2+(V_{(2)}\bm{q}^k)_j^2,\bm{x}_j\right) \\
& - \mu \left( (D_xV_{(1)}\bm{p}^k)_j^2  + (D_xV_{(2)}\bm{q}^k)_j^2 + (D_yV_{(1)}\bm{p}^k)_j^2 + (D_yV_{(2)}\bm{q}^k)_j^2 \right)].
\end{align*}

For the numerical tests, we consider the 2D NLS equation with $\mu = 1/2$ and
$$
R(|\psi|^2,x,y)= R_1(x,y)|\psi|^2 +  \frac{\beta}{2}|\psi|^4 , \quad R_1(x,y)=-\frac{1}{2}(x^2+y^2) - 2\exp (-(x^2+y^2)).
$$
In this form, the 2D NLS equation \eqref{2dnls1} is known as the  Gross-Pitaevskii equation.
For positive values of the parameter $\beta$, the solutions are focusing, and defocusing  for negative values.
We simulate the 2D NLS equation on the space-time domain $[-6,6]^2\times [0,10]$ by setting $\beta =1$.
The initial condition is calculated from the exact solution \cite{Li15nls}
$$
\psi(x,y,t) = \sqrt{2}\exp \left( -\frac{1}{2}(x^2+y^2) \right)\exp (-it).
$$
The simulations are performed on a space-time domain  with the spatial mesh sizes $\Delta x=\Delta y=12/100$ ($N=10000$), and the time-step size $\Delta t = 0.01$ ($N_t=1000$). With these setting, we obtain the snapshot matrices $\mathcal{P}$ and $\mathcal{Q}$ of size $10000\times 3000$.

\begin{figure}[H]
\centering
\includegraphics[width=0.5\columnwidth]{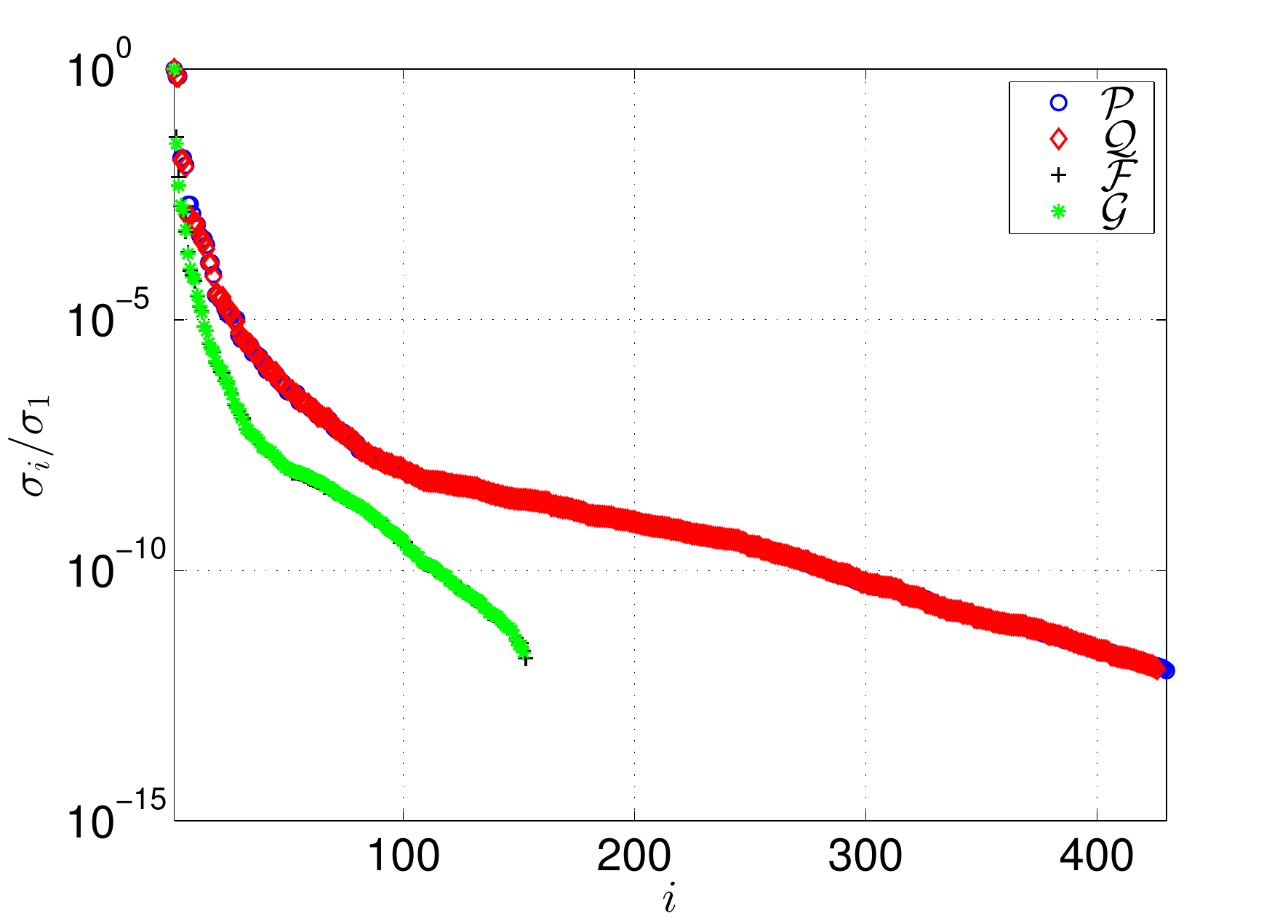}
\caption{2D NLS: Decay of singular values\label{nls2:sing}}
\end{figure}

The decay of the singular values of the solution snapshots $\mathcal{P}$, $\mathcal{Q}$ and the nonlinear snapshots $\mathcal{F}$, $\mathcal{G}$ is given in \figurename~\ref{nls2:sing}. The number of POD modes is taken as $n=10$, whereas $\tilde{n}=20$ for the DEIM modes.

\begin{figure}[H]
\centerline{\includegraphics[width=1.3\columnwidth]{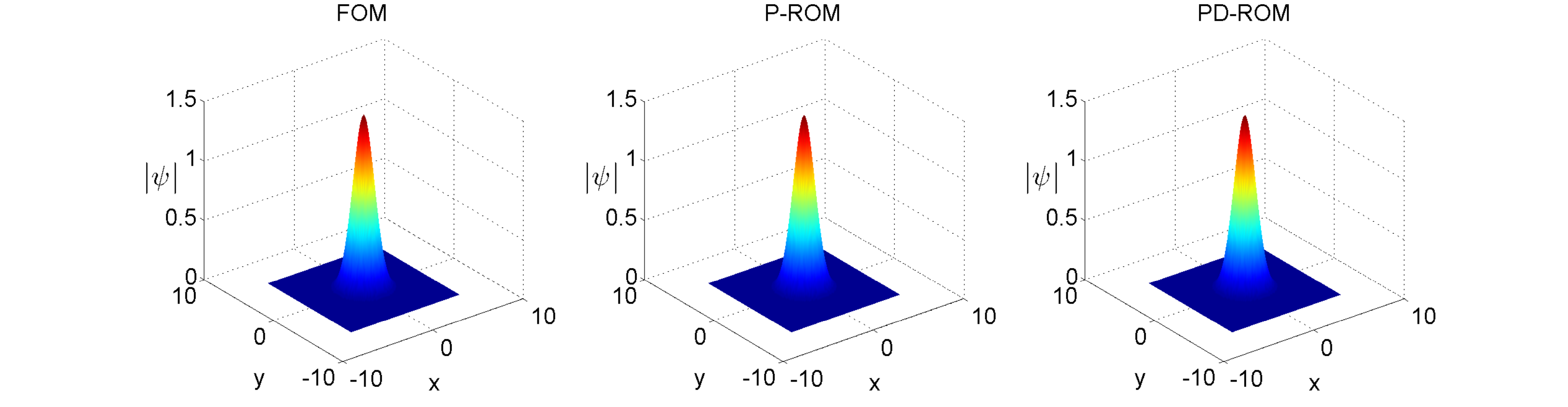}}
	\caption{2D NLS: Soliton waves at $t=5$\label{nls2:surf}}
\end{figure}

\begin{figure}[H]
\centering
\includegraphics[width=0.5\columnwidth]{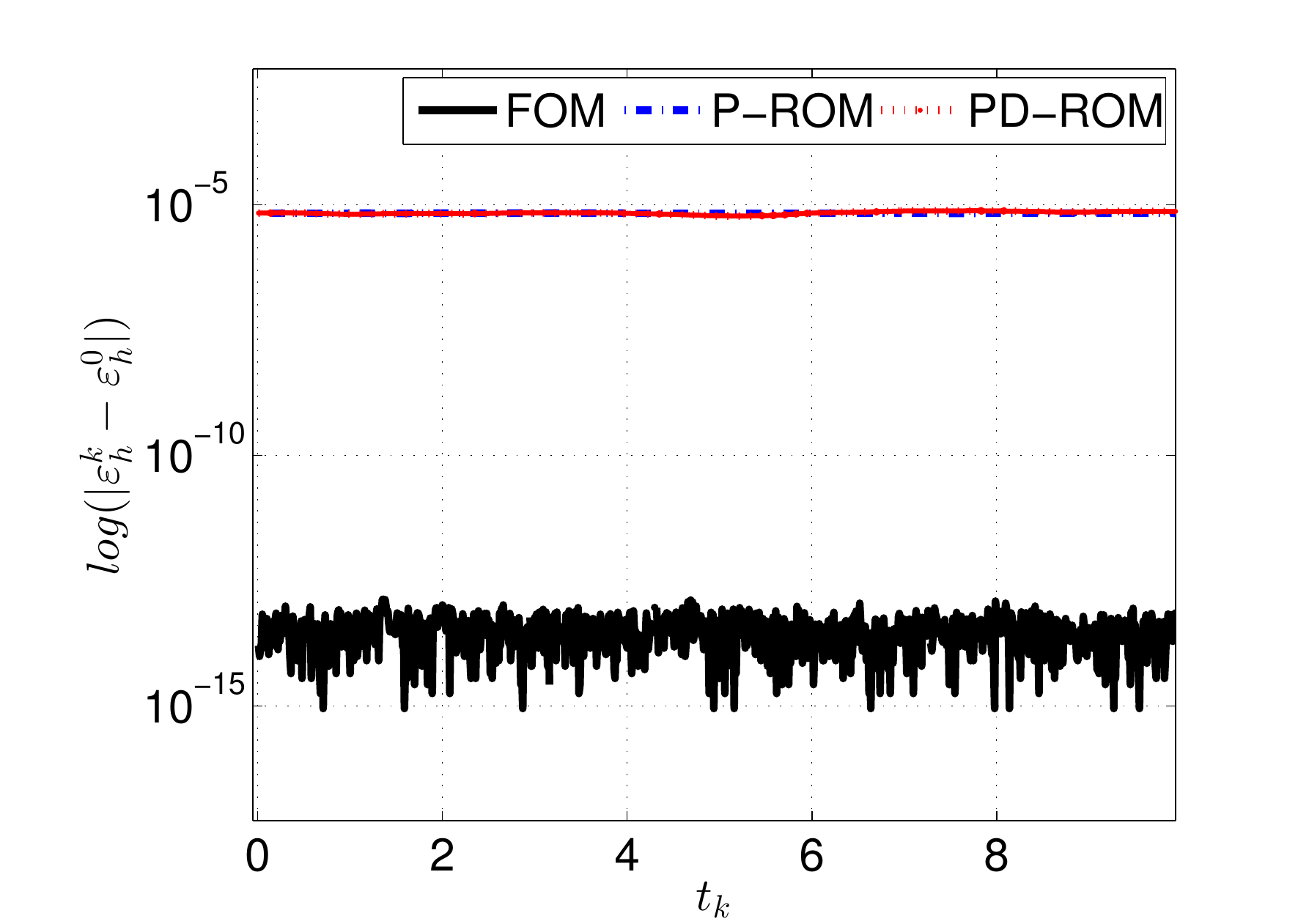}
\caption{2D NLS: Conservation of discrete global energy\label{nls2:energy}}
\end{figure}

The FOM solution and reduced approximations behave similar in \figurename~\ref{nls2:surf}.
The discrete reduced global energies are not with high accuracy as the full-order global energies in \figurename~\ref{nls2:energy} due to the reduced dimension of the ROMs.
When the reduced dimension is increased, the discrete global energy is preserved with higher accuracy, which incurs high computational cost, and the ROMs became computationally inefficient.
The oscillatory behavior of the discrete energies of the ROMs is an indication of the preservation of the global energy and stability of the solutions in long-term integration.

%%%%%%%%%%%%%%%%%%%%%%%%%%%%%%%%%%%%%%%%%%%%%%%%%%%%%%%%%%%%%%%%%%%%%%%%%%

Both ROMs have almost the same solution and shape errors as the FOM for the POD and DEIM modes determined with the prescribed tolerances in \tablename~\ref{tbl:error}.
The reduced global energies by the ROMs can not be preserved with the same accuracy as the FOM due to the reduced dimensions.
But they do not exhibit any drift over time, which reflects the energy-preserving nature of the reduced systems.

In \tablename~\ref{tbl:speed}, the computational efficiency of the ROMs is shown in terms of the wall-clock times in seconds and speedup factors of the reduced systems over the FOM.
The PD-ROMs have higher speedups than the P-ROMs, due to offline-online decomposition.
The computational gain by the ROMs is more pronounced in 2D examples due to large computing times of FOM.
Among the multi-symplectic equations, the NLS equation exhibit a lower speedup.
This results from the fact that in the NLS equation a coupled system with two unknown states have to be solved in contrast to the KdV equation.

\begin{table}[H]
  \centering
	\caption{Relative solution, relative shape and relative energy errors\label{tbl:error}}
\resizebox{\columnwidth}{!}{
\begin{tabular}{|l|ccc|ccc|ccc|}\hline
\multirow{2}{*}{Problem} & \multicolumn{3}{|c}{$ \text{E}_{\text{sol}} $} & \multicolumn{3}{|c}{$ \text{E}_{\text{shape}} $} & \multicolumn{3}{|c|}{$ \text{E}_{\text{energy}} $}  \\
 & FOM & P-ROM  & PD-ROM & FOM & P-ROM  & PD-ROM  & FOM & P-ROM  & PD-ROM \\
	\hline
KdV   &  4.82e-03    &  5.71e-03    & 5.62e-03    & 6.97e-05   & 9.67e-05  &  8.98e-05 & 1.78e-13 & 7.78e-06 & 9.90e-05  \\
	\hline
1D NLS   & 3.24e-02   & 3.42e-02   & 3.38e-02   & 9.70e-05 & 1.49e-04  &  1.75e-04   & 4.55e-15   & 1.27e-03   &  1.42e-03   \\
\hline
ZK   & 7.61e-03   &  7.64e-03  &  7.65e-03  & 8.42e-05  & 8.38e-05   & 8.39e-05    & 6.06e-14  &  2.65e-05  &  2.61e-05   \\
\hline
2D NLS   & 1.93e-02   & 1.92e-02   & 1.92e-02   & 4.83e-04  & 4.60e-04   &  4.60e-04   & 2.97e-14   &  1.47e-06  & 1.62e-06   \\
\hline
\end{tabular}}
\end{table}

\begin{table}[H]
\centering
\caption{Wall-clock times and speedup factors\label{tbl:speed}}
\resizebox{\columnwidth}{!}{
\begin{tabular}{|l|l|ccc|cc|}\hline
\multirow{2}{*}{Problem} & \multirow{2}{*}{\#Modes}  & \multicolumn{3}{|c}{Wall-Clock Time} & \multicolumn{2}{|c|}{Speedup}  \\
 & & FOM & P-ROM  & PD-ROM  & P-ROM  & PD-ROM \\
	\hline
KdV    & $ n=40 $, $\tilde{n}=45$ &  47.0  &  2.1  & 1.3  & 22.8  & 36.0 \\
	\hline
1D NLS & $ n=25 $, $\tilde{n}=45$ &  31.6  & 12.0  & 4.7  &  2.6  & 6.7 \\
\hline
ZK     & $ n=15 $, $\tilde{n}=25$ &  121.0 &  3.6  & 1.8   & 33.7  & 65.4 \\
\hline
2D NLS & $ n=10 $, $\tilde{n}=20$ & 273.9  & 22.8   & 6.8   &  12.0 & 40.3 \\
\hline
\end{tabular}}
\end{table}

In the following we compare our results with other structure-preserving ROMs for Hamiltonian systems. The FOMs are constructed by discretizing Hamiltonian PDEs by finite differences in space. Because in each paper, different time discretizations, initial conditions, and discretization parameters are used, an exact comparison is not possible.
Therefore, the comparison is carried out in terms of the reduced-order modes, preservation of the Hamiltonian, and computational efficiency, whenever provided.

\begin{itemize}

\item In \cite[Section 5.2.2]{Gong17}, the FOM and ROM have constructed for the one-dimensional KdV equation with a single soliton, and integrated in time by the AVF. When the reduced dimension increases from $n=40$ to $n=60$, the magnitude of the FOM-ROM error reduces from 5.64e-02 to 7.31e-04, and the associated maximum FOM-ROM Hamiltonian error decreases from 2.98e-04 to 9.15e-07.

In \cite[Section 5.1]{Miyatake19}, the KDV equation is integrated in time with the mid-point rule . The full-order solution profile is captured with $n=60$ POD modes. It was shown that the reduced Hamiltonian is well preserved even for small POD modes, for instance $n=20$.

In \cite[Section 4.1]{Karasozen21}, the FOM and ROM are obtained by a linearly implicit time integrator, i.e., Kahan's method. The full-order solutions are identical to reduced-order solutions for $n=30$ POD modes with the speedup factor $32$.

In \cite[Section 5.1.1]{Hesthaven21}, KdV equation with double soliton interaction is solved with a symplectic greedy ROM. As time integrators, the AVF and mid-point rules are used. The algorithm converges for a sufficiently large reduced space due to the complex wave propagation phenomena and slowly decaying Kolmogorov widths associated with the problem.

\item In the case of ZK equation, the ROM and FOM solutions are identical for $n=50$ modes, with a speedup factor $68$, whereas the full and reduced Hamiltonian errors are about order $10^{-3}$, as shown in \cite[Section 4.4]{Karasozen21}.

\item One-dimensional parametric NLS equation is integrated with the St{\"o}ormer-Verlet scheme in \cite[Section 5.2]{Hesthaven16}. The reduced system is obtained using the greedy algorithm with the cotangent lift, the complex SVD, the DEIM, the symplectic DEIM (SDEIM), and also the POD. The sizes of the ROMs are $n = 180$ for POD methods and $n = 90$ for symplectic methods. It is shown that the greedy, the cotangent lift, and the complex SVD methods generate stable ROMs that accurately approximate the full-order solutions.

\item The FOM of the two-dimensional NLS equation in parametric form is integrated in time with implicit mid-point rule, whereas the 2-stage partitioned RK method is used for the time integration of the ROM in \cite[Section 8.2]{Hesthaven22}. The adaptive dynamical ROM exhibits a speedup of $58$ over the FOM, with a smaller error than with the non-adaptive method.

\end{itemize}

Overall, the reduced-order solutions of the global energy preserving ROMs for multi-symplectic PDEs in \tablename~\ref{tbl:error}, have approximately the same level of FOM-ROM solution and energy errors as in the literature.  They are also computationally efficient as the structure-preserving ROMs.

%%%%%%%%%%%%%%%%%%%%%%%%%%%%%%%%%%%%%%%%%%%%%%%%%%%%%%%%%%%%%%%%%%%%%%%%%%%%%%%%%%%%%%%%%%%%
%%%%%%%%%%%%%%%%%%%%%%%%%%%%%%%%%%%%%%%%%%%%%%%%%%%%%%%%%%%%%%%%%%%%%%%%%%%%%%%%%%%%%%%%%%%%
\section{Conclusions}
\label{sec:conc}

In this paper, ROMs are constructed that preserve the global energy of well-known multi-symplectic PDEs.
We have proved that the discrete reduced global energy by the P-ROM is preserved exactly in the full-rank reduced space, and it is approximately preserved by the PD-ROM.
The reduced approximations are very close to the full-order solutions, and the shapes of the solitons are preserved by both ROMs.
The discrete reduced global energies are well preserved with small oscillations and do not show any drift in long-term integration.
All these reflect the global energy-preserving properties of the multi-symplectic ROMs.
Relatively large number of POD and DEIM modes show the limitation of the linear model order reduction techniques such as POD and DEIM for problems associated with transport and wave type phenomena as in this paper.
The linear reduced manifold can easily become large thus compromising ROMs efficiency.
Preservation of the local energy/momentum and multi-symplectic conservation in the ROM sense will be considered as the future research.

%%%%%%%%%%%%%%%%%%%%%%%%%%%%%%%%%%%%%%%%%%%%%%%%%%%%%%%%%%%%%%%%%%%%%%%%%%%%%%%%%%%%%%%%%%%%
%%%%%%%%%%%%%%%%%%%%%%%%%%%%%%%%%%%%%%%%%%%%%%%%%%%%%%%%%%%%%%%%%%%%%%%%%%%%%%%%%%%%%%%%%%%%
%\bibliographystyle{plain}
%\bibliography{references}

\end{document}